\documentclass[aos,numbers,noinfoline]{imsart}
\usepackage{amsmath}
\usepackage{amsthm,hyperref}
\usepackage{amssymb,amsfonts}
\usepackage[mathscr]{eucal}
\usepackage{graphicx}
\usepackage{color}
\usepackage{enumerate}
\usepackage{bm}
\usepackage[comma,sort&compress]{natbib}
\usepackage{multirow}

\usepackage{booktabs}
\usepackage{bigstrut}

\newtheorem{theorem}{Theorem}
\newtheorem{lemma}{Lemma}

\numberwithin{theorem}{section}
\numberwithin{lemma}{section}
\numberwithin{corollary}{section}

\newcommand{\mZ}{\mathbb{Z}}

\def\ci{\perp\!\!\!\perp}

\newcommand{\ex}{\ensuremath{{\mathbb E}}}

\newcommand{\E}{\ensuremath{{\mathbb E}}}

\newcommand{\X}{\ensuremath{\bm{X}}}
\newcommand{\Y}{\ensuremath{\bm{Y}}}
\newcommand{\x}{\ensuremath{\bm{x}}}
\newcommand{\y}{\ensuremath{\bm{y}}}
\newcommand{\thet}{\ensuremath \bm{\theta}}

\newcommand{\phat}{\ensuremath{\hat{p}}}

\newcommand{\thetab}{\ensuremath \bm{\thetab}}

\newcommand{\RR}{\mathbb{R}}

\newcommand{\hf}{\textstyle{\frac{1}{2}}}

\newcommand{\g}{\textsf{G}}
\newcommand{\s}{\textsf{S}}
\newcommand{\bg}{\textsf{B}}

\usepackage{xr}
\begin{document}

\begin{frontmatter}
	
	\title{On Minimax Optimality of Sparse Bayes Predictive Density Estimates}
	\runtitle{Sparse Bayes Predictive Density Estimates}

	\begin{aug}
\author{\fnms{Gourab} \snm{Mukherjee}\ead[label=e1]{gourab@usc.edu}}
\address{\printead{e1}}
\and
\author{\fnms{Iain M.} \snm{Johnstone}\ead[label=e2]{imj@stanford.edu}}
\address{\printead{e2}}
\affiliation{University of Southern California and Stanford University}
\hspace{1cm}\\
\end{aug}

\begin{abstract}
We study predictive density estimation under
Kullback-Leibler loss in $\ell_0$-sparse Gaussian sequence models. We
propose proper Bayes predictive density estimates and establish
asymptotic minimaxity in sparse models. 

A surprise is the existence of a phase transition in the
future-to-past variance ratio $r$. For $r < r_0 = (\surd 5 - 1)/4$,
the natural discrete prior ceases to be asymptotically optimal.
Instead, for subcritical $r$, a `bi-grid' prior with a central region
of reduced grid spacing recovers asymptotic minimaxity.
This phenomenon seems to have no analog in the otherwise 
parallel  theory of point estimation of a multivariate normal mean
under quadratic loss.

For spike-and-slab priors to have any prospect of minimaxity, we show
that the sparse parameter space needs also to be magnitude constrained.  
Within a substantial range of magnitudes, spike-and-slab
priors can attain asymptotic minimaxity.
\end{abstract}
	

	
\end{frontmatter}

\section{Introduction and Main results}\label{sec:intro}

Predictive density estimation is a
fundamental problem in statistical prediction analysis
\citep{Aitchison-book,Geisser-book}.  
Here, it is studied in a 
high dimensional Gaussian setting under sparsity assumptions on the
unknown location parameters. 
Fuller references and background for the problem are 
given after a formulation of our main results.

 We consider a simple Gaussian model for high dimensional prediction:
\begin{equation}
\label{eq:modelM2}
\X \sim N_n(\thet, v_x \bm{I}), \qquad
\Y \sim N_n(\thet, v_y \bm{I}), \qquad
\X \ci \Y | \thet.
\end{equation}
Our goal is to predict the distribution of a future observation $\Y$
on the basis of the `past' observation vector $\X$. In this model, the
past and future observations are independent, but are linked by the
common mean parameter $\thet$ which is assumed to be unknown.  The
variances $v_x$ and $v_y$ may differ and are assumed to be known. 

The true probability densities of $\X$ and $\Y$ are denoted by
$p(\x| \thet, v_x)$ and $p(\y| \thet, v_y)$ respectively. We seek
estimators $\hat p(\y|\x)$ of the future observation density
$p(\y|\thet, v_y)$, and study their risk properties under sparsity
assumptions on $\thet$ as dimension $n$ increases to $\infty$.


To evaluate the performance of a predictive density estimator $\hat
p(\y|\x)$, we use Kullback-Leibler `distance' as loss
function:
\begin{displaymath}
L(\thet, \hat p(\cdot|\x))
= \int p(\y|\thet , v_y) \log \frac{p(\y|\thet , v_y)}{\hat p(\y|\x)} d\y.
\end{displaymath}
The corresponding KL risk function follows by averaging over the
distribution of the past observation:
\begin{displaymath}
\rho(\thet, \hat p)
= \int L(\thet, \hat p(\cdot|\x) p(\x|\thet , v_x)) d \x.
\end{displaymath}
Now, given a prior measure $\pi( d\thet)$, the average or integrated risk is
\begin{equation}
\label{eq:integ-risk}
B(\pi, \hat p)
= \int \rho(\thet, \hat p) \pi(d \thet).
\end{equation}
For any prior measure $\pi( d \thet)$, proper or improper,
such that the posterior $\pi(d \thet | \x)$ is well defined, the Bayes
predictive density is given by
\begin{equation}
\label{eq:bpd-temp}
\hat p_\pi(\y|\x)
= \int p(\y|\thet, v_y) \pi(d \thet|\x).
\end{equation}
The Bayes predictive density in (\ref{eq:bpd-temp}) minimizes both the
posterior expected loss
$\int L( \thet, \hat p(\cdot |\x)) \pi(d \thet |\x)$ and the
integrated risk $B(\pi, \hat p)$ in the class of all density estimates,
e.g. Sec. 2.4 of \citep{Mukherjee-thesis}.
The minimum is
the Bayes KL risk:
\begin{align}\label{def-bayes-risk}
B(\pi) := \inf_{\hat p} B(\pi, \hat p)~.
\end{align}
\par
We study the predictive risk $\rho( \thet,\hat p)$ in a high
dimensional setting under an $\ell_0$-sparsity condition on the
parameter space. This `exact' sparsity condition has been widely used
in statistical estimation problems, e.g. \citep[Ch. 8]{Johnstone-book}. With
$\| \thet \|_0 = \# \{ i: \theta_i \neq 0 \}$,  consider
the parameter set:
\begin{displaymath}
\Theta_n[s] = \{ \thet \in \mathbb{R}^n: \| \thet \|_0 \leq s \}.
\end{displaymath}
The minimax KL risk for estimation over $\Theta$ is given by
\begin{equation}
\label{eq:mmx-risk-def}
R_N(\Theta) = \inf_{\hat p} \sup_{\thet \in \Theta} \ \rho( \thet,
\hat p),
\end{equation}
the infimum being taken over \textit{all} predictive density
estimators $\hat p(\y|\x)$. We often write $\mathsf{pde}$ for 
predictive density estimate. The notation $a_n \sim b_n$ denotes
$a_n/b_n \to 1$ as $n \to \infty$ and $a_n=O(b_n)$ denotes $|a_n/b_n|$ is bounded for all large $n$.

\subsection{Main Results.} 
Henceforth, we assume $v_x =1$. As the problem is scale equivariant,
results for general $v_x$ will easily follow. 
A key parameter is the future-to-past variance ratio
\begin{equation}
\label{eq:oracle}
r = v_y/v_x = v_y, \qquad   v = (1+r^{-1})^{-1}.
\end{equation}
Here $v$ is 
the `oracle variance'
which would be the variance of the UMVUE for $\thet$, if both $\X$ and
$\Y$ were observed. The variance ratio $r$
determines not only the magnitude of the
minimax risk but also the construction of minimax optimal
\textsf{pdes.}
In our asymptotic model, the dimension $n \to \infty$ and the
sparsity $s = s_n$ may depend on $n$, but the variance ratio $r$
remains fixed.  

In the sparse limit
$\eta_n = s_n/n \to 0$, for any fixed $r \in (0,\infty)$,
\citet{Mukherjee-15} evaluated the minimax risk to be:
\begin{equation}\label{eq:minimax-risk}
R_N(\Theta_n[s_n])\sim \frac{1}{1+r} s_n \log (n/s_n) = \frac{1}{1+r} n \eta_n \log \eta_n^{-1},
\end{equation}
and a thresholding based \textsf{pde} was shown to
attain the minimax risk. 

By their nature, thresholding rules are not smooth functions of the
data.
This paper develops proper Bayes \textsf{pdes}
-- necessarily smooth functions --
that are asymptotically minimax in sparse regimes.
Throughout we consider sparse symmetric priors
\begin{equation}
   \label{eq:ssp}
  \pi[\eta] = (1-\eta) \delta_0 + \hf \eta(\nu^+ + \nu^-),
\end{equation}
where $\delta_0$ is unit mass at $0$, and 
$\eta \in [0,1]$ is the sparsity parameter,
while $\nu^+$ is a probability
measure on $(0,\infty)$ and $\nu^-$ is its reflection on
$(-\infty,0)$. 
Priors on $\thet$ are built from i.i.d. draws:
\begin{equation*}
  \pi_n (d \thet) = \prod_{i=1}^n \pi[\eta_n] (d \theta_i),
\end{equation*}
where $\eta_n=s_n/n$.
The Bayes \textsf{pde} based on prior $\pi_n$ is the product
density estimate: 
$$\phat(\y|\x)=\prod_{i=1}^n \phat(y_i|x_i)~.$$

\begin{figure}[b] 
	\includegraphics[width=0.95\textwidth]{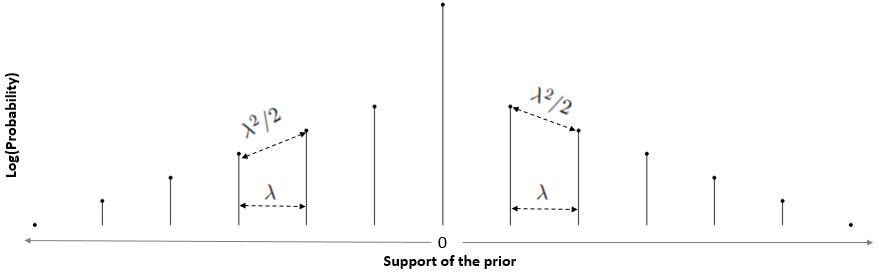}
	\caption{\small \sl Schematic for the grid
          prior. The uniform spacing $\lambda$ between the support
          points is shown on 
          the x-axis and the logarithm of the probabilities of the
          support points on the y-axis, the latter having linear decay.
} 
\label{fig-0a}
\end{figure}

We begin with a discrete `grid prior' $\nu_\g^+$ in which the support
points have equal spacing
\begin{equation}
    \label{eq:lambdadefs}
  \lambda = \lambda(\eta) = \sqrt{2 \log \eta^{-v}},
\end{equation}
and geometric mass decay at rate $\eta^v = e^{-\lambda^2/2}$.
More precisely,
\begin{equation*}
  \nu_\g^+ = (1-\eta^v) \sum_{j=1}^\infty \eta^{(j-1)v}
  \delta_{\lambda j}.
\end{equation*}
The corresponding sparse grid prior $\pi_\g[\eta]$ built via \eqref{eq:ssp}
has a schematic illustration in Figure \ref{fig-0a}. 
Such discrete priors are a natural starting point for our predictive
setting given their optimality properties in point estimation,
recalled in the next subsection.

Our first result gives a precise description of the first order
asymptotic maximum risk of the Bayes \textsf{pde} $\hat{p}_\g$ based on the multivariate product prior $\pi_{\g,n}$. 
Define 
\begin{equation}\label{eq:define.hr}
  \begin{split}
     h_r & =(1+2r)(1+r)^{-2}(1-2r-4r^2)/4  \ \leq \ 1/4 \\
     h_r^+ & =\max(h_r,0).
  \end{split}
\end{equation}
Let $r_0=(\sqrt{5}-1)/4$ be the positive root of the equation
$4r^2+2r-1=0$, and note that $h_r^+ > 0$ iff $r <r_0$. 

\begin{theorem}\label{th:D-minimax}
	As $\eta_n = s_n/n \to 0$, for any fixed $r \in (0,\infty)$ we have
	\begin{displaymath}
	\max_{\Theta_n[s_n]} \rho(\thet, \hat p_{\g}) = R_N(\Theta_n[s_n])\big(1+h_r^+ + o(1)\big) \text{ as } n \to \infty~.
	\end{displaymath}
\end{theorem}
Thus for all $r \geq r_0$, $\phat_{\g}$ is exactly minimax
optimal, while for all $r < r_0$, it is minimax suboptimal but still attains the
minimax rate, and has maximum risk at most 1.25 times the minimax
value, whatever be the value of $r$. 

As the future-to-past variance ratio $r$ decreases,
the difficulty of the
predictive density estimation problem increases,
as we have to estimate the future observation density based
on increasingly noisy past observations.
Theorem \ref{th:D-minimax} shows that
rules which are minimax optimal for higher values of $r$ can be
sub-optimal for lower values of $r$.
This phenomenon was seen with threshold density estimates 
in \citep[Sec. S.2, Lemma~S.2.1]{Mukherjee-15} 
as well as in the recent work of
\citep{maruyama2016harmonic} on non-sparse prediction.

\begin{figure}[b] 
	\includegraphics[width=0.95\textwidth]{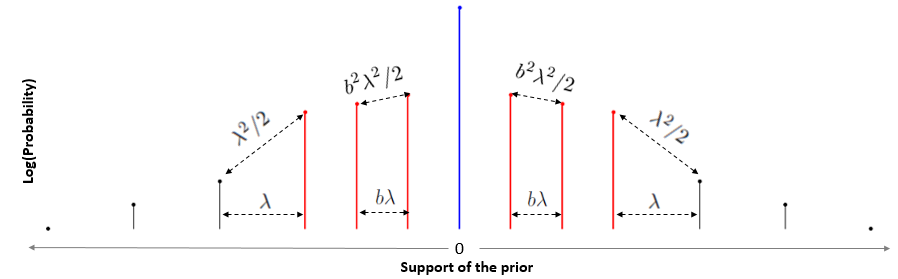}
	\caption{\small \sl Schematic for the bi-grid
          prior. The $x$-axis now shows the two spacings, and the 
$y$-axis the two different rates of log-linear decay of the prior
probabilities. 
} \label{fig-0b}
\end{figure}

To obtain asymptotic minimaxity for all $r$, we need to modify the prior.
The \textsf{Bi-grid} $\pi_{\bg}$ prior is obtained from $\pi_G$ by selecting an `inner zone'
on which the spacing of the prior atoms is reduced from $\lambda$ to
$b \lambda$, where
\begin{equation*}
  b = \min \{ 4r(1+r)(1+2r)^{-1},1\}.
\end{equation*}
Note that $b < 1$ iff $r < r_0$.
The decay ratio in the inner zone is increased from $\eta^v =
e^{-\lambda^2/2}$ to $\eta^{vb^2} = e^{-b^2 \lambda^2/2}$. 
See Figure~\ref{fig-0b} for a schematic depiction.

More precisely, $\pi_B$ is a sparse symmetric prior
of form \eqref{eq:ssp} with
\begin{equation*}
  \nu_B^+
   = c \Big[ \sum_{k=1}^K \eta^{(k-1)vb^2} \delta_{\nu_k} 
              + \eta^{(K-1)vb^2} \sum_{j=1}^\infty \eta^{jv}
            \delta_{\mu_j} \Big].
\end{equation*}
The normalization $c = c(\eta)$ is at \eqref{eq:normconst}.
The support points fall in two zones:
\begin{itemize}
\item[(i)] Inner zone: \qquad $\nu_k = \lambda + (k-1)b\lambda$ for $k
  = 1, \ldots, K$ 
\item[(ii)] Outer zone: \qquad $\mu_j = \nu_K + j\lambda$ for $j = 1,
  2, \ldots$ 
\end{itemize}
The cardinality of the inner zone is
\begin{equation}
    \label{eq:def.K}
  K = 1 + \lceil 2b^{-3/2} \rceil.
\end{equation}

A main result of the paper is that the Bayes predictive density
estimate $\phat_{\bg}$ based on the product prior $\pi_{\bg,n}$ is 
asymptotically minimax optimal.  
\begin{theorem}\label{th:M-minimax}
For each fixed $r \in (0,\infty)$, as $\eta_n = s_n/n \to 0$,  we have
	\begin{displaymath}
	\max_{\Theta_n[s_n]} \rho(\thet, \hat p_{\bg}) =
        R_N(\Theta_n[s_n])(1+o(1)) \qquad \text{ as } n \to \infty.
	\end{displaymath}
\end{theorem}
The following theorem shows that the bi-grid prior $\pi_{\bg}$ is also
asymptotically least favorable. 
\begin{theorem}\label{th:M-least-fav}
	If $s_n \to \infty$ and $s_n/n \to 0$, then
	\begin{displaymath}
	B(\pi_{\bg,n}) = R_N(\Theta_n[s_n]) \cdot (1 + o(1)).
	\end{displaymath}
\end{theorem}

Unlike Theorem~\ref{th:M-minimax} we need the assumption that
$s_n \to \infty$. It ensures that $\pi_{\bg,n}$
actually concentrates on $\Theta_n[s_n]$, namely that
$\pi_{\bg,n}(\Theta_n[s_n]) \to 1$ as $n \to \infty$.  For the case
where $s_n$ does not diverge to $\infty$ an asymptotically least
favorable prior can be constructed from a sparse prior built from
`independent blocks'. The construction is discussed in
Section~\ref{sec-2-last-pf}.

\subsection{Discussion}
\label{sec:discussion}
A fully Bayesian approach is a natural route to \textsf{pdes} with
good properties \citep{Hartigan98,Aslan06}, with advantages over 
`plug-in' or thresholding based density estimates.
Indeed, a coordinatewise threshold rule 
$\phat_{\text{T}}(\y|\x)=\prod_{i=1}^n \phat_{\text{T}}(y_i|x_i)$
is typically built from univariate \textsf{pdes} which combine two Bayes
\textsf{pdes} --  for example based on 
uniform $\phat_{\text{U}}$ and cluster priors $\phat_{\text{CL}}$, as
in   
\citep[Eq. (14)]{Mukherjee-15}:
\begin{equation*}
\phat_{\text{T}}(y_i|x_i)=
\begin{cases}
\phat_{\text{U}}(y_i|x_i) & \text{if } |x_i| > v^{-1/2}\lambda \\
\phat_{\text{CL}}(y_i|x_i) & \text{if } |x_i| \leq  v^{-1/2} \lambda. 
\end{cases}
\end{equation*}
This is manifestly discontinuous as a function of the data $\x$.

The bi-grid Bayes rule achieves the same purposes as the hybrid
$\phat_{\text{T}}$. Indeed, the close spacing $b \lambda$ in the inner
section of $\pi_\bg$ yields the same risk control as the (unevenly
spaced) cluster prior for small and moderate $\theta$, while the
uniform $\lambda$ spacing in the outer section of $\pi_\bg$ controls
risk for large $\theta$ in the same way as the uniform prior.

Decision theoretic parallels between predictive density estimation and
the point estimation of a Gaussian mean under quadratic loss have been
established by
\citep{George06,George08,George12,Brown08,Komaki01,Kobayashi08,Xu12,Ghosh08}
for unconstrained $\thet$, and by 
\citep{Xu10}, \citep{Fourdrinier11},
\citep{kubokawa2013minimaxity} and \citep{Mukherjee-15} for various
constraint sets $\Theta$.

The phase transition seen in Theorems \ref{th:D-minimax} and
\ref{th:M-minimax} seems however to have no parallel in point
estimation.  Indeed, it follows from \citep{Johnstone94a} that a first
order minimax rule for quadratic loss in the sparse setting is derived
from the Mallows prior \citep{mall78}, with
$\nu_{\mathsf{Q}}^+ = (1-\eta)\sum_{j=1}^\infty \eta^{j-1}
\delta_{\lambda_e j}$.  Here
$\lambda_e = \sqrt{2 \log \eta^{-1}} = v^{-1/2}\lambda$ so that the
predictive setting involves a reduced spacing in the prior.
More significantly, there is no analog in point estimation of the
inner section with its 
further reduced spacing for $r < r_0$. 


Our main technical contribution lies in sharp methods for bounding the
global KL risk for general bi-grid priors, see Lemmas \ref{lem:bounds}
and \ref{lem.pf.D3}, and for spike-and-slab priors, Section
\ref{sec-2-3}. 
The sharp predictive
risk bounds established here provide new asymptotic perspectives in
the information geometric framework of
\cite{Komaki96,komaki2006shrinkage,yano2017information} and augment
new sparse prediction techniques for general multivariate
predictive density estimation theory developed in
\cite{George06,Brown08,Komaki01,kubokawa2015predictive,kubokawa2017predictive,matsuda2015singular}.

\subsection{ Minimax risk of Spike and Slab priors}
\label{sec:minimax-risk-spike}

Some of the most popular Bayesian variable selection techniques
are built on the ``spike and slab" priors  
\citep{mitchell1988bayesian,george1997approaches,ishwaran2005spike}.
Such priors and their computationally tractable extensions have found
success in variable selection in high-dimensional sparse regression
models, e.g.
\citep{park2008bayesian,bhattacharya2015dirichlet,rockova2014spike,rockova2015bayesian,rovckova2014negotiating,ishwaran2005spike-a}
and the references therein. 
While this is a well established methodological research area
\citep{o2009review},
optimality of their respective predictive density estimates
has so far not been studied.

Here, we consider simple ``spike and slab" prior distributions in the
flavor of the foundational paper \citep{mitchell1988bayesian}.
Thus, let
\begin{equation}\label{eq:def-SS-prior}
\pi_{\s}[\eta,l] = (1-\eta) \delta_0 +  \eta/(2l) I\{\mu \in [-l,l]\}
\, \, d\mu ~. 
\end{equation}
Now consider the multivariate Bayes predictive density $\phat_S[l]$
based on the $n$ dimensional prior containing i.i.d. copies of
$\pi_S[s_n/n,l]$ as before. 
Such Bayes \textsf{pdes} are necessarily asymptotically sub-optimal.
\begin{lemma}\label{lem-uncont}
For any fixed $l \in [0,\infty)$, we have
$$\lim_{n \to \infty} \bigg\{\max_{\Theta_n[s_n]} \rho(\thet, \hat p_S[l]) \bigg \}\bigg/ R_N(\Theta_n[s_n]) = \infty~ $$
\end{lemma}
The result is hardly surprising, as the support of $\pi_\s$ is restricted to 
$[-l,l]$ for $l$ fixed, and the corresponding \textsf{pde} has large
risk away from the support.
Consider therefore bounded subsets of the sparse parameter sets
$\Theta_n[s_n]$:
\begin{displaymath}
\Theta_n[s,t] = \{ \thet \in \mathbb{R}^n:  \| \thet \|_0 \leq s \text{ and } |\theta_i| \leq t \text{ for all } i =1,\ldots,n\}.
\end{displaymath}

We allow  $t = t_n$ to increase with $n$, and note next that the
increase must be at least as fast as $\lambda_n$ to have minimax risk
equivalent to $\Theta_n[s_n]$.

\begin{lemma}
  \label{lem:tn-risk}
For all $t_n$ there is a simple bound
\begin{equation*}
  R_N(\Theta_n[s_n,t_n]) \leq s_n t_n^2/(2r).
\end{equation*}
If $t_n > \lambda_n = \sqrt{2 \log \eta_n^{-v}}$, then
\begin{equation}
  \label{eq:asy-equiv}
  R_N(\Theta_n[s_n,t_n])
   \sim s_n \lambda_n^2/(2r)
   \sim R_N(\Theta_n[s_n]).
\end{equation}
\end{lemma}

The next result says that there is a substantial range of sparsities 
$\eta_n = s_n/n$ for which $\hat p_\s[t_n]$ is asymptotically minimax
over the bounded sets $\Theta_n[s_n,t_n]$, and hence also over
$\Theta_n[s_n]$:
\begin{theorem}\label{thm-cont-1}
	As $\eta_n=s_n/n \to 0$,
suppose that $t_n/(\log \eta_n^{-1})^{1/2} \to
        \infty$ but $\log t_n/(\log \eta_n^{-1}) \to 0$.
Then
	$$\lim_{n \to \infty} \bigg\{\max_{\Theta_n[s_n,t_n]}
        \rho(\thet, \hat p_\s[t_n]) \bigg\} \bigg /
        R_N(\Theta_n[s_n,t_n]) = 1~.$$ 
\end{theorem}

Note that if  $t_n \to \infty$ at a rate slower than 
$(\log \eta_n^{-1})^{1/2}$ then, by Lemma~\ref{lem:tn-risk}, 
$R_N(\Theta_n[s_n,t_n])$ is no longer 
equivalent to $R_N(\Theta_n[s_n])$ as $n \to \infty$.
At the other extreme,  we show next that if $t_n$ grows at rate 
$\eta_n^{-\beta}$ or higher for any $\beta > 0$, then 
no spike and uniform slab procedure can be 
minimax optimal. 
\begin{theorem}\label{thm-cont-2}
If $\eta_n = s_n/n \to 0$ and $\log t_n = \beta
\log \eta_n^{-1}$ for some $\beta > 0$, then as $n \to \infty$,
\begin{equation*}
\min_{l>1} \max_{\Theta[s_n,t_n]}
\rho(\bm{\theta},\hat{p}_S[l])
\geq (1 + \beta)  R_N(\Theta_n[s_n,t_n]) (1 + o(1))~.
\end{equation*}
\end{theorem}

\subsection{Organization of the Paper} 
Section \ref{sec-2} presents the risk properties of the Grid
and Bi-grid prior based \textsf{pdes} and proofs of the main
results. Section~\ref{sec-2-3} proves the
spike-and-slab results. Section~\ref{sec:simu} compares the
performance of the \textsf{pdes} through simulation experiments.
Proofs of all lemmas are in
the Appendix.


\section{Proof of the main results}\label{sec-2}

We focus on priors with i.i.d. components
$\pi(d \bm{\theta}) = \prod_{i=1}^n \pi_1(d \theta_i)$.
The predictive density then has product form
$\hat{p}_\pi (\mathbf{y}|\mathbf{x}) = \prod_{i=1}^n
\hat{p}_{\pi_1}(y_i|x_i)$ 
and the predictive risk is additive
\begin{equation}
\label{eq:additive}
\rho(\bm{\theta},\hat{p}_\pi)
= \sum_{i=1}^n \rho(\theta_i,\hat{p}_{\pi_1}).
\end{equation}
When the context is clear, we drop the subscript and write
$\hat{p}_\pi$ for the univariate Bayes predictive density also.

For our sparse parameter sets $\Theta_n[s]$ and $\Theta_n[s,t]$, there
is an easy reduction of the maximum multivariate risk of a product
rule to a univariate risk maximum.
Indeed, \eqref{eq:additive} yields 
\begin{equation}
\label{eq:mult}
s_n \sup_{|\theta_n|\leq t_n} \rho(\theta, \hat{p}_1)
\leq \sup_{\Theta_n[s_n,t_n]} \rho(\bm{\theta},\hat{p})
\leq n(1-\eta_n) \rho(0,\hat{p}_1) 
+ s_n \sup_{\theta \in \RR}  \rho(\theta, \hat{p}_1).
\end{equation}

\textit{Sparse priors.} \ 
We turn now to the predictive risk properties of univariate priors of
the form
\begin{equation*}
\pi(d\mu) = (1-\eta) \delta_0 + \eta \nu(d \mu).  
\end{equation*}
The following risk decomposition is fundamental; it is proved in the appendix.

\begin{lemma} \label{lem:decomp}
	Let $Z \sim \mathcal{N}(0,1)$.
	For a sparse prior,
	\begin{equation}  \label{eq:dec}
	\rho(\theta, \hat{p}_\pi ) 
	= \frac{\theta^2}{2r} - \E \log N_{\theta,v}(Z) + \E \log D_\theta(Z),
	\end{equation}
	where $D_\theta(Z) = N_{\theta,1}(Z)$ and
	\begin{equation}
	\label{eq:Nthetav}
	N_{\theta,v}(Z) 
	= 1 + \frac{\eta}{1-\eta} \int \exp \left\{ \frac{\mu Z}{\sqrt v} 
	+ \frac{\mu \theta}{v} -
	\frac{\mu^2}{2v} \right \} \nu(d \mu)~.
	\end{equation}
\end{lemma}

\bigskip
Clearly $N_{\theta,v}(Z), D_\theta(Z) > 1$, and so we have the simple but useful 
``basic lower'' and ``basic upper'' risk bounds 
\begin{equation}
\label{eq:triv}
\frac{\theta^2}{2r} - \E \log N_{\theta,v}(Z) 
\leq \rho(\theta, \hat{p}_\pi ) 
\leq \frac{\theta^2}{2r} + \E \log D_\theta(Z).
\end{equation}

From Jensen's inequality,
\begin{equation}
\label{eq:jensen}
\E \log N_{\theta,v}(Z)  \leq  \log \big(\E N_{\theta,v}(Z)\big),
\end{equation}
and since $\E \exp (\zeta Z)= \exp (\zeta^2/2)$,
\begin{equation}
\label{eq:mean}
\E N_{\theta,v}(Z)
=  1 + \frac{\eta}{1-\eta} \int \exp \Big( \frac{\mu \theta}{v} 
\Big) \nu(d \mu),
\end{equation}
and, in particular,
\begin{equation*}
\E D_0(Z) = \E N_{0,v}(Z) = (1-\eta)^{-1}.
\end{equation*}

Consequently, from the right side of \eqref{eq:triv}, then
\eqref{eq:jensen} (for $v=1$) and 
the previous display,
\begin{equation}   \label{eq:sparse-at-0}
\rho(0,\hat{p}_\pi) 
\leq \log (1-\eta)^{-1} 
= \eta (1 + o(1)) \quad \text{as} \quad \eta \to 0.
\end{equation}

Return now to the univariate reduction \eqref{eq:mult}.
From \eqref{eq:sparse-at-0} it is clear that
$n\rho(0,\hat{p}_1) \leq n \eta_n(1+o(1)) = s_n(1+o(1))$.
So for the minimaxity results of Theorems ~\ref{th:D-minimax}, ~\ref{th:M-minimax} and ~\ref{thm-cont-1}, it
suffices to show the univariate bound 
\begin{equation}
\label{eq:univ-bd}
\sup_{\theta \in \RR}  \rho(\theta, \hat{p}_1)
\leq \lambda_n^2/(2r) +o(\lambda_n^2),
\end{equation}
for then
\begin{equation*}
\sup_{\Theta_n} \rho(\bm{\theta}, \hat{p}_\pi) \leq 
s_n [\lambda_n^2/(2r) + o(\lambda_n^2)].
\end{equation*}

\subsection{Risk Properties of the Grid and the Bi-grid Priors based estimates}\label{sec:2-proofs} 
To allow a unified analysis, we introduce a class of discrete sparse
priors that includes both grid and bi-grid priors.
For $0 < b \leq 1$ and $r > 0$, let 
\begin{equation}
\label{eq:prior}
\pi_\mathsf{D}[\eta;b,r] = \sum_{j \in \mZ} \pi_j \delta_{\mu_j}
\end{equation}
where
$\mu_{-j} = - \mu_j$, $\pi_{-j} = \pi_j$.
The support points satisfy $\mu_0 = 0$  and $\mu_j = \lambda \alpha_j$ for $j > 0$,
where the piecewise linear spacing function
\begin{equation}\label{eq:def.spacing}
\alpha_j =
\begin{cases}
1 + b(j-1)  & \quad 1 \leq j \leq K \\
\alpha_K + j-K & \quad \quad \ \ j > K
\end{cases}
\end{equation}
has increments $\dot \alpha_j = \alpha_{j+1} - \alpha_j = b$ or $1$
according as $j \leq K$ or $j>K$. 
Set $\zeta=\eta^{v}$.
The prior masses are given by
\begin{equation*}
\pi_0 = 1- \eta, \qquad
\pi_j = c(\eta) \eta \zeta^{\beta_j -1}, 
\end{equation*}
for $j \geq 1$. 
The decay function in the prior probabilities 
\begin{equation}\label{eq:def.decay}
\beta_j =
\begin{cases}
1 + b^2(j-1)  & \quad 1 \leq j \leq K \\
\beta_K + j-K & \quad \quad \ \ j > K
\end{cases}  
\end{equation}
has the same form as $\alpha_j$ with $b$ replaced by $b^2$.
This choice is crucial for Lemma~\ref{lem:bounds} below and its
consequent risk bounds.
In particular, note that 
$\beta_j \leq \alpha_j$ and that the
increments $\dot \beta_j = \beta_{j+1} - \beta_j$
satisfy
\begin{equation}
\label{eq:scaling}
\dot \beta_j  = \dot \alpha_j^2 \qquad \quad \text{all  } j \geq 1.  
\end{equation}
In addition, $l \to \alpha_l^2 - \beta_l$ is increasing for $l \geq
1$, as
\begin{equation}
\label{eq:increase}
(d/dl)(\alpha_l^2 - \beta_l) = \dot \alpha_l(2\alpha_l - \dot
\alpha_l) \geq 0.
\end{equation}
The normalizing constant $c(\eta)$ is determined by
\begin{equation}
\label{eq:normconst}
\frac{1}{2c(\eta)}
= \frac{1 - \eta^{b^2 vK}}{1 - \eta^{b^2 v}}
+\frac{\eta^{b^2 v(K-1)+v}}{1-\eta^v}.
\end{equation}

Since $\pi_\mathsf{D}$ is a sparse prior, we may apply the
decomposition of predictive risk given in Lemma \ref{lem:decomp}. 
Inserting the discrete measure \eqref{eq:prior}, we obtain
\begin{align}
N_{\theta,v}(Z) 
& = 1 + \sum_{j \neq 0} N_j, \label{eq:N}\\
N_j & = {\pi}_0^{-1}\,{\pi_j}\, \exp \{ v^{-1/2} \mu_j Z +
v^{-1} (\mu_j \theta - \hf \mu_j^2 ) \} \label{eq:Nj}
\intertext{In the special case $v=1$, it will be helpful to write
	$D_\theta(Z) = N_{\theta,1}(Z)$ as} 
D_\theta(Z) & = 1 + \sum_{j \neq 0} D_j,   \label{eq:Dtheta}   \\
D_j & = {\pi}_0^{-1}\,{\pi_j}\, \exp \{ \mu_j Z + \mu_j \theta -
\hf \mu_j^2  \}.
\label{eq:Aj}
\end{align}

The probability ratio $\pi_j/\pi_0$ can also be written in exponential
form. 
To this end, introduce $c_1(\eta) = c(\eta) (1-\eta)^{-1}$.
Recall that $v^{-1} = 1 + r^{-1}$ and $\zeta = \eta^v = \exp
(-\lambda^2/2)$ and then
rewrite $\eta = \zeta^{v^{-1}} = \exp \{ -\hf \lambda^2(1+r^{-1})
\}$. We arrive at
\begin{equation}
\label{eq:piratio}
{\pi}_0^{-1}\,{\pi_j} 
= c_1(\eta) \exp \{ - \hf \lambda^2(\beta_j + r^{-1}) \}.
\end{equation}
We can therefore, for example, rewrite
\begin{equation}
\label{eq:altDrep}
\begin{split}
D_j & = \exp \{ \mu_j Z - G(\mu_j;\theta) \} \\
G(\mu_j;\theta) & = \hf \mu_j^2 - \mu_j \theta + \hf \lambda^2
(\beta_j + r^{-1}).
\end{split}
\end{equation}

To obtain an upper bound for $\rho(\theta,\hat{p}_\mathsf{D})$ we use
\eqref{eq:Nthetav}. 
We focus on two consecutive terms $N_{j}$, $N_{j+1}$ in \eqref{eq:N};
ignoring all other terms trivially yields a lower bound for
$N_{\theta,v}$.
For the upper bound for $D_\theta$, 
a single (suitably chosen) term $D_j$ in \eqref{eq:Dtheta} suffices,
but more care is needed 
to show that the neglected terms are negligible. 


Bring in a co-ordinate system $(l, \omega)$ for $\theta$: each
$\theta \geq 0$ 
can be uniquely written in the form $\theta = \lambda(\alpha_l +
\omega)$ for a uniquely determined $l \in \mathbf{N}$ and $\omega \in
[0,\dot \alpha_l)$.
We can therefore write $l = l(\theta)$ and $\omega = \omega(\theta)$.

We argue heuristically that $l(\theta)$ is an appropriate choice of
index for our bounds.
Indeed, from \eqref{eq:Nj} and \eqref{eq:piratio},
\begin{equation}
\label{eq:Elog}
\E \log N_j = c - \hf \{ (\mu_j-\theta)^2/v -
\lambda^2 \beta_j \}
\end{equation}
after collecting terms not involving $j$ into $c$. 
Hence, for $\theta \in [\mu_l,\mu_{l+1})$, the choice $j = l$ or $l+1$ will
minimize or nearly minimize the quadratic,
and these suffice for the lower bound.
For $D_\theta$, we have from \eqref{eq:altDrep} that
$\E \log D_j = - G(\mu_j; \theta)$.
We show in the Appendix (in the proof of Lemma~\ref{lem:bounds})
that $j \to G(\mu_j;\theta)$ is indeed minimized at $j = l$
for \textit{each} $\theta \in [\mu_l,\mu_{l+1})$.

Focus therefore on the terms $N_{l(\theta)}$ and $D_{l(\theta)}$. 
When $\theta = \lambda(\alpha_l + \omega)$, 
\begin{equation*}
\mu_j \theta - \hf \mu_j^2 
= \hf \lambda^2 ( 2 \alpha_j(\alpha_l+ \omega) - \alpha_j^2).
\end{equation*}
Combining this with \eqref{eq:piratio}, 
for $j = l, l+1$, we can write
\begin{equation}
\label{eq:Al}
\begin{split}
N_l & = c_1(\eta) \exp \{ 
\hf \lambda^2 n(l,\omega) +  \alpha_l \lambda Z /\sqrt{v} \}  \\
N_{l+1} & = c_1(\eta) \exp \{ 
\hf \lambda^2 \check n(l,\omega) +  \alpha_{l+1} \lambda Z /\sqrt{v} \}  \\
D_l & = c_1(\eta) \exp \{          
\hf \lambda^2 d(l,\omega) +  \alpha_l \lambda Z \}      
\end{split}
\end{equation}
in terms of three linear functions of $\omega$:
\begin{equation}\label{eq:n}
\begin{split}
n(l,\omega) & = v^{-1}(\alpha_l^2 + 2 \alpha_l \omega) - \beta_l -
r^{-1} \\ 
d(l,\omega) & = \qquad \alpha_l^2 + 2 \alpha_l \omega - \beta_l -
r^{-1}.
\end{split}
\end{equation}
and, corresponding to $N_{l+1}$,
\begin{equation}
  \label{eq:ndiff}
  \check{n}(l,\omega)
     = n(l,\omega) + 2v^{-1}\dot \alpha_l \omega - (1+v^{-1})\dot \alpha_l^2.
\end{equation}

We now state our key uniform bounds on the risk components of
\eqref{eq:prior}. 

\begin{lemma}
\label{lem:bounds}
For any fixed $r \in (0,\infty)$ and $b \in [0,1]$,
with $\lambda$ defined in \eqref{eq:lambdadefs},
uniformly in $\theta = \lambda(\alpha_l + \omega) \geq \lambda$, we
have the following bounds: 
\begin{align*} 
\E \log N_{\theta,v}(Z) 
& \geq \hf \lambda^2 (n \vee \check{n})(l,\omega) + O(1), \\
\E \log D_{\theta}(Z)
& \leq \hf \lambda^2 d^+(l,\omega) + O(\lambda), 
\end{align*}
For $0 \leq \theta < \lambda$ we just have
$\E \log N_{\theta,v}(Z) \geq 0$, and 
$\E \log D_{\theta}(Z) \leq O(\lambda).$
\end{lemma}

The rather intricate proof is given in the appendix. 
The appearance of the positive part of $d(l,\omega)$ in the upper
bound may be understood this way:
if $d(l,\omega) < 0$, we cannot expect the term $D_l$ to dominate $D_0
= 1$ in \eqref{eq:Dtheta}.

In the reverse direction, we need only a bound for $\theta \in 
[\mu_1,\mu_2]$ in our proofs of theorems~\ref{th:D-minimax} and \ref{th:M-minimax}.  
\begin{lemma}
  \label{lem.pf.D3}
For any fixed $r \in (0,\infty)$ and $b \in [0,1]$,
with $\lambda$ defined in \eqref{eq:lambdadefs},
uniformly in $\theta \in \lambda[\alpha_1,\alpha_2]$, we have
\begin{equation*}
  \E \log N_{\theta,v}(Z)
    \leq \hf \lambda^2 n(1,\omega) + O(\lambda).
\end{equation*}
\end{lemma}

\subsection{Proof of Theorems~\ref{th:D-minimax} and \ref{th:M-minimax}}
Inserting the bounds of Lemma \ref{lem:bounds} in risk decomposition
\eqref{eq:dec}, we get
\begin{align}
\rho(\theta,\hat{p}_\mathsf{D})
& = (2r)^{-1} \lambda^2 \sigma(l,\omega)  + O(\lambda) \notag \\
\sigma(l,\omega)
& =  (\alpha_l+\omega)^2 - r (n\vee \check n)(l,\omega) + r
d^+(l,\omega).
\label{eq:rhobd}  
\end{align}

Our task is to investigate when $\sigma(l,\omega) \leq 1$. 
First assume $d(l,\omega) \geq 0$. In this easier case,
 from \eqref{eq:n}, $r(n-d) = \alpha_l^2 + 2 \alpha_l \omega$,
and so
\begin{equation*}
  \sigma(l,\omega) \leq (\alpha_l+\omega)^2 - r(n-d) 
                   = \omega^2 \leq 1,
\end{equation*}

Now, suppose that $d(l,\omega) < 0$.
For $l =0$, control on the risk is immediate from \eqref{eq:triv},
and so, from now on $l \geq 1$.  
We compare $n(l,\omega)$ and $\check
n(l,\omega)$ for $\omega \in [0,\dot \alpha_l]$
by using \eqref{eq:ndiff}.
Both are linear functions of $\omega$,
intersecting at $\omega_* = \dot \alpha_l(1+v)/2 < \dot \alpha_l$. 
Now 
$n(l,0) > \check n(l,0)$ while $\check n$ has a larger positive slope. 
Hence $n \vee \check n$ equals $n$  on $[0,\omega_*]$ and
$\check n$ on $[\omega_*,\dot \alpha_l]$.
Observe that on the right interval $[\omega_*,\dot \alpha_l]$,
\begin{equation*}
  \partial \sigma/\partial \omega (l,\omega) 
    =2[\alpha_l + \omega - r v^{-1} (\dot \alpha_l + \alpha_l)]
    \leq 2(1-rv^{-1})\alpha_{l+1}
    < 0.
\end{equation*}
Hence, in seeking the maximum of $\sigma(l,\omega)$ on $[0,\dot
\alpha_l]$, we may confine attention to $0 \leq \omega \leq
\omega_*$.  On this range
\begin{equation}
   \label{eq:sig-bd}
   \begin{split}
  \sigma(l,\omega) 
     & = (\alpha_l+\omega)^2 - r n(l,\omega) \\
     & = \omega^2 + (\alpha_l^2+2 \alpha_l \omega )(1-rv^{-1}) + r
       \beta_r + 1 \\
     & = 1 + \omega^2 +r(\beta_l-\alpha_l^2) - 2r \alpha_l \omega  \\
     & \leq 1 + \omega^2 - 2r\omega = \bar \sigma(\omega),     
   \end{split}
\end{equation}
say, where we used $\alpha_l \geq 1$ and 
$\alpha_l^2 -\beta_l \geq \alpha_1^2 - \beta_1 = 0$, from
\eqref{eq:increase}.

In particular $\bar \sigma(0) = 1$ and $\bar \sigma(\omega_*) = 1 +
\omega_*(\omega_*-2r).$ 

For the grid prior, $b=1$ and one evaluates 
\begin{equation*}
  \bar \sigma(\omega_*) = 1 + (1+2r)(1+r)^{-2}(1-2r-4r^2)/4= 1 + h_r~.
\end{equation*}
Consequently, 
$\rho(\theta,\hat{p}_\mathsf{G}) \leq (2r)^{-1} \lambda^2 \{1+h_r^+\}  +
	O(\lambda)~,$
which establishes the upper bound in Theorem \ref{th:D-minimax}.
For the lower bound, we look at the risk at $\theta_*=\lambda(1+\omega_*)$. 
Apply Lemma \ref{lem.pf.D3} using 
$n(1,\omega_*) = 2 v^{-1} \omega_*$, to get
from \eqref{eq:triv}
\begin{equation*}
  \rho(\theta_*,\hat{p}_\mathsf{G}) 
    \geq (2r)^{-1}\lambda^2 \{(1+\omega_*)^2-2r v^{-1} \omega_* \} +O(\lambda)
    = (2r)^{-1} \lambda^2 (1+h_r^+) +O(\lambda),
\end{equation*}
since the quantity in braces equals $1+\omega_*^2-2r\omega_* = \bar
\sigma(\omega_*) = 1 + h_r^+.$ 
This completes the proof of Theorem \ref{th:D-minimax}.

We now turn our attention to proving Theorem~\ref{th:M-minimax}. 
We first verify that if $b \leq \min \{1, 4r\}$, then
$l \geq K$ necessarily implies $d(l,\omega) \geq 0$. 
Set $K_* = 1 + 2 \,b^{-3/2}.$
From the monotonicity \eqref{eq:increase}, along with 
$\omega \geq 0$, we have
\begin{align*}
d(l,\omega) 
& \geq \alpha_l^2 - \beta_l - r^{-1}
\geq  \alpha_{K_*}^2 - \beta_{K_*} - r^{-1} \\
&=  4 b^{-1/2} -  2 b^{1/2} + 4 b^{-1}  -r^{-1}	\geq 0,
\end{align*}
using $b \leq 1$ for the first two terms and $b \leq 4r$ for the
second pair.

Now return to \eqref{eq:sig-bd} -- we have
\begin{equation*}
  \sigma(l,\omega) \leq 1 + \omega_*(\omega_*-2r)_+.
\end{equation*}
Since $d(l,\omega) \leq 0$, we must
have $l < K$ and hence $\dot \alpha_l=b$. 
Clearly $\sigma(l,\omega) \leq 1$ so long as 
$\omega_* = b(1+v)/2 \leq 2r$, or equivalently, 
$b \leq 4r/(1+v)$.
So in this case, for \textit{all} $\theta$ we have 
$\rho(\theta,\hat{p}_\mathsf{D}) \leq (2r)^{-1} \lambda^2 +
O(\lambda)$, which establishes \eqref{eq:univ-bd} and hence Theorem~\ref{th:M-minimax}.

\subsection{Proof of Theorem~\ref{th:M-least-fav}}\label{sec-2-last-pf}
As $\pi_{\bg,n}$ is i.i.d. and due to the product structure of the problem, its Bayes risk simplifies
$$B(\pi_{\bg,n}, \phat_\bg) = n B(\pi_\bg,\phat_\bg).$$
For the univariate problem the Bayes risk of the prior $\pi_{\bg}$ is
\begin{align*}
	B(\pi_\bg, \phat_\bg) &\geq  c_0 \eta_n \, \{\rho(\lambda_n,\phat_\bg) + \rho(-\lambda_n,\phat_\bg)\}\\
	&=2 c_0 \eta_n \, \rho(\lambda_n,\phat_\bg)
	 \geq   2 c_0 \eta_n \, [\lambda_n^2/(2r) - \ex \log
          N_{\lambda_n,v}(Z) ], 
\end{align*}
where the equality above follows by symmetry and the inequality by
\eqref{eq:dec}. 
From 
\eqref{eq:normconst} we have $2 c_0 \geq 1 - O(\eta_n^{b^2 v})$. 
Lemma~\ref{lem.pf.D3} shows that 
$\ex\log N_{\lambda,v}(Z) = O(\lambda)$ because $n(1,0)=0$. 
Hence $B(\pi_\bg, \phat_\bg) \geq \eta_n \lambda_n^2/(2r) \cdot
(1+o(1))$ and the proof is done.

 
\textit{Remark.} When $s_n$ does not diverge to $\infty$,  an
`independent blocks' sparse prior using $\pi_{\bg}$ is asymptotically
least favorable, along the lines of [\citealp{Johnstone-book}, Ch. 8.6].
Let $\pi_S(\tau ; m)$ denote a single spike prior of scale $\tau$ on
$\mathbb{R}^m$. This chooses an index $I \in \{1,\ldots,m\}$ at random and
sets $\theta=\tau e_I$, where $e_i$ is a unit length vector in the
$i$th co-ordinate direction. 
We randomly draw $\tau$ from
$(\nu_{\bg}^+ + \nu_{\bg}^-)/2$. 
However, instead of
\eqref{eq:lambdadefs}, we choose $\lambda = v^{1/2}(t_m - \log t_m)$
where $t_m=\sqrt{2 \log m}$.  
The independent blocks prior $\pi_{\sf IB,n}$ on $\Theta[s_n]$ is built by dividing $\{1,\ldots,n\}$
into $s_n$ contiguous blocks $B_j$, each of length $m = m_n =
[n/s_n]$.
Independently for each block $B_j$, draw components according to
$\pi_S(\cdot; m)$ and set $\theta_i = 0$ for the remaining $n -
m_ns_n$ coordinates. This prior 
is supported on $\Theta[s_n]$ as any draw from $\pi_{\sf IB,n}$ has exactly $s_n$ non-zero
components. The proof that it is least favorable is then analogous 
to that of Theorem 6 in \citep{Mukherjee-15}.


\section{Risk properties of Spike and Slab procedures}\label{sec-2-3}
We again use the risk decomposition provided by
Lemma~\ref{lem:decomp}, now with the univariate spike and slab prior
$\pi_S[\eta,l]$. We use $N^S_{\theta,v}(Z)$ and $D^S_{\theta}(Z)$ to
denote the associated risk components of Lemma~\ref{lem:decomp} for
the spike and slab predictive density estimates $\phat_S[l]$ based on
the prior $\pi_S[\eta,l]$ for some $l >0$ (the dependence on $l$ is
kept implicit in the notations).

\smallskip
\textit{Proof of Lemma \ref{lem:tn-risk}.} \ 
For the first upper bound, simply take $\nu = \delta_0$ in Lemma
\ref{lem:decomp}; the corresponding $\pi_0 = \delta_0$ has 
$\rho(\theta,\hat{p}_{\pi_0}) = \theta^2/(2r)$. The bound now follows
from \eqref{eq:mult}. 
For the second statement, we claim that whenever $t_n > \lambda_n$, then as
$n \to \infty$, 
\begin{equation}
  \label{eq:equiv}
  R_N(\Theta_n[s_n,t_n])
   \sim R_N(\Theta_n[s_n])
   \sim s_n \lambda_n^2/(2r).
\end{equation}
Indeed, the independent blocks prior $\pi_n^{IB}$ constructed in
\citep[Theorem 6]{Mukherjee-15}
to show that $R_N(\Theta_n[s_n]) \sim s_n \lambda_n^2/(2r)$
is actually, by its very definition, supported on
$\Theta_n[s_n,\nu_n]$,
where $\nu_n < \sqrt{v} \sqrt{2 \log [n/s_n]} \leq \lambda_n < t_n$.
Since obviously 
$\Theta_n[s_n,\nu_n] \subset \Theta_n[s_n,t_n] \subset\Theta_n[s_n]$, 
the conclusion (\ref{eq:equiv}) follows. \qed.

\smallskip
For lower bounds on risk of its predictive density estimate, the following convexity inequality is helpful. It
is proved in the appendix. 

\begin{lemma}\label{lem:sec3-1}
If $\eta \leq \hf$ and $\theta l / v \geq 1$, then
\begin{equation*}
\E \log N^S_{\theta,v}(Z) \leq \theta l/v.
\end{equation*}
\end{lemma}  
The \textit{proof of Lemma~\ref{lem-uncont}} follows easily from the above lemma.  
From the left side of (\ref{eq:triv}) and Lemma 3.1,
\begin{equation*}
  \rho(\theta, \hat{p}_S[l]) \geq \frac{\theta^2}{2r} - \frac{\theta
    l}{v}, \qquad \quad \text{for } \quad \theta \geq \frac{v}{l}.
\end{equation*}
Hence,
\begin{equation*}
  \sup_{\Theta_n[s_n]} \rho(\mathbf{\bm{\theta}}, \hat{p}_S[l])
       = s_n \sup_{\theta \in \RR} \rho(\theta, \hat{p}_S[l]) 
       = \infty.
\end{equation*}

\subsection{Proof of Theorem~\ref{thm-cont-1}}
An alternative representation for $N_{\theta,v}^S$ will be useful.
Completing the square in \eqref{eq:Nthetav}, we get
\begin{equation}
  \label{eq:Nrep}
  N_{\theta,v}^S(Z) = 1 + c(\eta) \sqrt{v} \exp (\hf Z_{\theta,v}^2)
  \Phi_{l,v}, 
\end{equation}
where we have set $Z_{\theta,v} = Z + \theta/\sqrt{v}$ and
\begin{equation*}
  \Phi_{l,v} = \Phi(v^{-1/2}(l-\theta)-Z) - 
               \Phi(v^{-1/2}(-l-\theta)-Z).
\end{equation*}
In the appendix, we show that, uniformly in $v \in (0,1)$, $l \geq 1$
and $|\theta| \leq l$,
\begin{equation}
  \label{eq:Phi-bd}
  \E \log \Phi_{l,v} \geq a_0 := \log \phi(0) + 2/3.
\end{equation}

The constant $c(\eta) = \eta(1-\eta)^{-1} \{2l \phi(0) \}^{-1}$
satisfies
\begin{equation}
  \label{eq:ceta-bds}
  - \log l - \lambda^2/(2v) \leq 
    \log \{(1-\eta)^{-1} c(\eta) \}
    \leq   \hf \log \tfrac{\pi}{2} - \log l - \lambda^2/(2v)
\end{equation}

From the preceding three displays and $\E Z_{\theta,v}^2 = 1 +
\theta^2/v$ we obtain
\begin{align}
  - \E \log N_{\theta,v}^S (Z)
  & \leq - \log c(\eta) - \hf \log v - \hf \E Z_{\theta,v}^2
          - \E \log \Phi_{l,v} \notag \\
  & \leq \log l + \lambda^2/(2v) - \theta^2/(2v) + O(1). \label{eq:bd1}
\end{align}

For $D^S_\theta$, we show in the appendix that for each $r> 0$, and
with $\lambda = \sqrt{2 v \log
  \eta^{-1}}$ and $\tilde{\lambda} = \lambda/\sqrt{v} + \sqrt{2 \log
  \lambda}$,
\begin{equation}
  \label{eq:D-bd}
  \E \log D_\theta^S(Z)
    \leq
    \begin{cases}
      O(\lambda \log \lambda)      & 0 < \theta < \tilde{\lambda} \\
      \theta^2/2 - \lambda^2/(2v) + O(\lambda) & \theta \geq
      \tilde{\lambda}.
    \end{cases}
\end{equation}

We assemble these pieces to bound $\rho(\theta,\hat{p}_S[l])$ and claim that
\begin{equation*}
  \rho(\theta,\hat{p}_S[l]) \leq
  \begin{cases}
    \theta^2/(2r) + O(\lambda \log \lambda) & 0 < \theta < \lambda \\
    \lambda^2/(2r) + \log l + O(\lambda \log \lambda) & 
            \lambda \leq \theta < \tilde{\lambda} \\
    \log l + O(\lambda)    & \tilde{\lambda} \leq \theta.
  \end{cases}
\end{equation*}
For $0 < \theta < \lambda$, simply use the basic upper bound
(\ref{eq:triv}) along with (\ref{eq:D-bd}).
For the remaining two cases, we use the full
decomposition (\ref{eq:dec}) of Lemma \ref{lem:decomp}. 
To this end, note from \eqref{eq:bd1} and
$v^{-1} = r^{-1}+1$ that
\begin{equation*}
  \theta^2/(2r) - \ex \log N_{\theta,v}^S(Z)
   \leq \lambda^2/(2r) - (\theta^2-\lambda^2)/2 + \log l + O(1).
\end{equation*}
Combining this with the bounds in \eqref{eq:D-bd} yields the
remaining two bounds.

For any $l \geq 1$ such that $\log l = o(\lambda^2)$, we conclude that
as $\lambda \to \infty$,
\begin{equation*}
  \sup_\theta \rho(\theta, \hat{p}_S[l]) 
       \leq \frac{\lambda^2}{2r} (1 + o(1)).
\end{equation*}
This completes the proof of (\ref{eq:univ-bd}) and, as remarked there,
the proof of Theorem \ref{thm-cont-1}.

\subsection{Proof of Theorem~\ref{thm-cont-2}}
We use the basic lower risk bound \eqref{eq:triv}, and show that for
suitable $\theta$ that $\E \log N_{\theta,v}^S{Z}$ cannot be large
enough to offset the leading term $\theta^2/(2r)$. 
To obtain a result uniform over all slab widths $l$, we need two
different types of upper bound on $N_{\theta,v}^S$.

Define $t_\lambda$ and $\tilde{t}_\lambda = o(t_\lambda)$ by setting
$\log t_\lambda = \beta \lambda^2/(2v)$
and $\log \tilde{t}_\lambda = \log t_\lambda - \lambda.$
We look first at large values of $l$, using representation
(\ref{eq:Nrep}).
Observe first that for $l > \tilde{t}_\lambda$, the right side of
(\ref{eq:ceta-bds}) yields
\begin{equation*}
  \sqrt{v} c(\eta) 
    \leq C \exp \{ - \log \tilde{t}_\lambda - \lambda^2/(2v) \}
    = C \exp \{ - \tilde{\theta}^2/(2v) \}
\end{equation*}
for a constant $C = C(v)$ if we set 
$\tilde{\theta}^2 = \lambda^2 + 2v \log \tilde{t}_\lambda$. 
Using now (\ref{eq:Nrep}) and $\Phi_{l,v} < 1$, we have
\begin{align*}
  \log N_{\tilde{\theta},v}^S(Z)
   & \leq \log \{ 1 + C \exp [-\tilde{\theta}^2/(2v) +
     (Z+\tilde{\theta}/\sqrt{v})^2/2] \} \\
   & \leq \log 2 + \log(1+C) + Z^2/2 + |Z| \tilde{\theta}/\sqrt{v}.
\end{align*}
Consequently $\E \log N_{\tilde{\theta},v}^S(Z)\leq k_1 + k_2
\tilde{\theta}$ where $k_i = k_i(v)$. 
Hence, from the left side of risk bound (\ref{eq:triv}),
\begin{equation*}
  \rho(\tilde{\theta}, \hat{p}_S[l]) 
    \geq \frac{\tilde{\theta}^2}{2r} - k_1 \tilde{\theta} - k_2.
\end{equation*}
Now observe  from the definition of $\tilde{t}_\lambda$ that
$\tilde{\theta}^2 = (1+\beta)\lambda^2 - 2v \lambda$ and that
$\tilde{\theta} < t_\lambda$ for large $\lambda$.
We conclude that for large $\lambda$,
\begin{equation}
  \label{eq:small-l}
  \inf_{l > \tilde{t}_\lambda} \sup_{\theta \in [0,t_\lambda]}
    \frac{\rho(\theta,\hat{p}_S[l])}{\lambda^2/(2r)} 
    \geq 1 + \beta + O(\lambda^{-1}).
\end{equation}

For $l \leq \tilde{t}_\lambda$, we set $\theta = t_\lambda$ and use
the left side of (\ref{eq:triv}), then Lemma 3.1: 
\begin{equation*}
  \sup_{\theta \leq t_\lambda} \rho(t_\lambda,p_S[l]) 
      \geq \frac{t_\lambda^2}{2r} - \frac{t_\lambda l}{v}
      \geq \frac{t_\lambda^2}{2r} - \frac{t_\lambda
        \tilde{t}_\lambda}{v}
      \geq \frac{t_\lambda^2}{2r} (1 + o(1)),
\end{equation*}
where in the last inequality we used $\tilde{t}_\lambda =
o(t_\lambda)$.
Consequently,
\begin{equation}
  \label{eq:large-l}
  \inf_{l \leq \tilde{t}_\lambda} \sup_{\theta \in [0,t_\lambda]}
    \frac{\rho(\theta,\hat{p}_S[l])}{\lambda^2/(2r)} 
    \geq \frac{r t_\lambda^2}{\lambda^2} (1 + o(1)).
\end{equation}
Combining (\ref{eq:small-l}) with (\ref{eq:large-l}) and then using 
(\ref{eq:mult}) to go over to the multivariate problem, we obtain
\begin{equation*}
  \min_{l>1} \max_{\Theta_n[s_n,t_n]} \rho(\bm{\theta}, \hat{p}_S[l]) 
    \geq (1 + \beta) s_n \lambda_n^2/(2r) (1 + o(1)).
\end{equation*}
Theorem~\ref{thm-cont-2} now follows from \eqref{eq:asy-equiv} of 
Lemma \ref{lem:tn-risk}.



\section{Numerical Experiments}\label{sec:simu}


We looked  at the numerical effectiveness of our
asymptotic results under different levels of sparsity $\eta_n$, with
special focus on moderate values. 
The product structure  and the good bounds
\eqref{eq:mult} relating maximal multivariate and univariate risks
allow us to concentrate on the maximal risk of the 
univariate \textsf{pdes}. We use a constrained prior space
\begin{equation*}
\mathfrak{m}_l(\eta)=\{\pi \in \mathcal{P}(\RR): \pi(\theta = 0) \ge
1-\eta,\, \pi(|\theta| > l)=0\},
\end{equation*}
and set $l = 5 \lambda = 5 \sqrt{2 \log \eta^{-v}}$.
We consider three sparsity levels:
(a) Moderate: $\eta = 0.1$, (b) High: $\eta =
0.001$, (c) Very High: $\eta = 10^{-10}$.

\begin{figure}[!h] 
	\includegraphics[width=\textwidth]{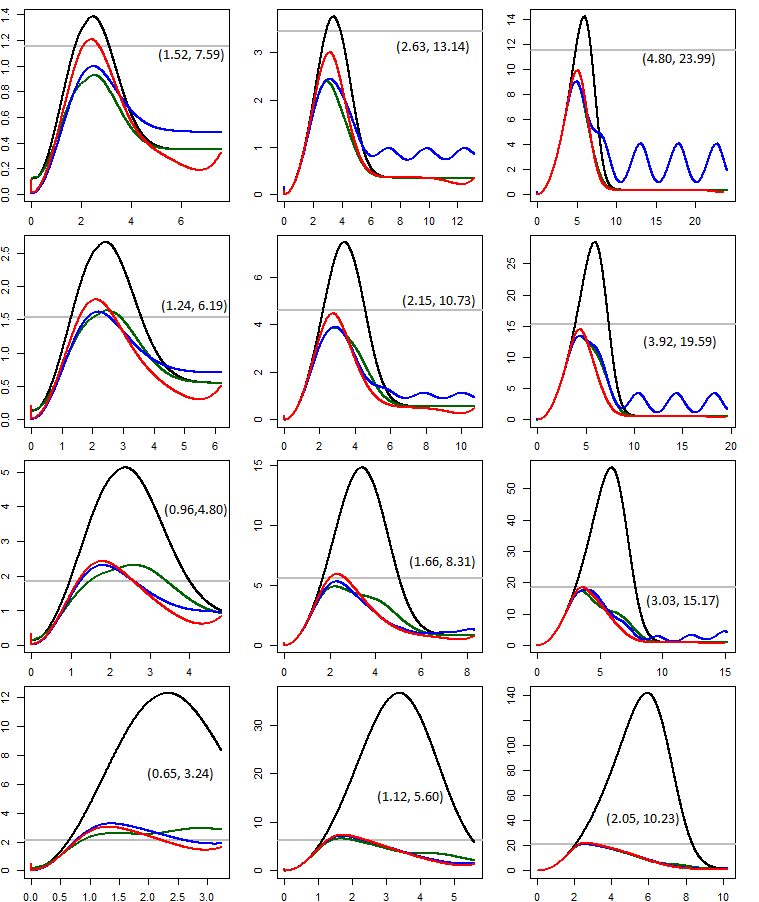}
\caption{\small \sl 
Risk plots $\rho(\theta, \cdot )$ for univariate predictive density
estimators  $\phat_H$ (black), $\phat_T$ (dark green), $\phat_{\bg}$
(blue) and $\phat_S$ (red) versus $\theta \in [0, l]$, for $l = 5 \lambda$.
Columns vary with moderate, high and very high sparsity,
$\eta=0.1, 0.001, 10^{-10}$, left to right. 
Rows vary $r=1, 0.5, 0.25$ and $0.1$ from top to bottom.
The horizontal line shows the asymptotic univariate minimax
          risk of $\log \eta^{-1}/(1+r) = \lambda^2/(2r)$,
with $\lambda = \sqrt{2 \log \eta^{-v}}$ and $l$ shown in the insets.}
\label{fig.compare.risk}
\end{figure}

Figure~\ref{fig.compare.risk} shows univariate risk plots in the three 
sparsity regimes 
for the following \textsf{pdes}:
\begin{itemize}
	\item Hard threshold Plug-in  \textsf{pde} (H-Plugin):
	$$\phat_H(y|x)=p(y|\hat \theta_H,v_y) \text{ where }  \hat
        \theta_H(x) = x \, I\{|x|>(v_x/v)^{1/2}\lambda\}~.$$ 
	\item Cluster prior and Thresholding (C-Thresh) based
          asymptotically minimax \textsf{pde} $\phat_T$ proposed in
          \cite[Eqn. (12)-(14)]{Mukherjee-15}  
	\item Bayes \textsf{pde} based on the bi-grid prior (Bi-Grid): $\phat_{\bg}$ 
	\item Spike and Slab  predictive density estimator (SS): $\phat_S[l]$ .
\end{itemize}
The basic features of the risk plots are unchanged
even under moderate sparsity. The hard threshold plug-in density
estimator $\phat_H$ does poorly for small values of $r$. 
For each $r$, the
maximal risks of $\phat_T$ and $\phat_{\bg}$ lie near or below the
asymptotic level of $\log \eta^{-1}/(1+r)$  under
high and very high sparsity, and at worst moderately above the
asymptotic level for moderate sparsity.

Table~\ref{numerical.risk.1} reports the maximum value of the risk
plots for these predictive estimators.  [Table~\ref{numerical.risk.2}
in the Supplement shows the locations of their respective maxima.]
The tables and plots show that the Bi-grid prior Bayes
\textsf{pde} and the C-Thresh
\textsf{pde} $\phat_T$ have similar worst case performance. 
The
maximal risk of the spike and slab procedure is higher than that of
$\phat_T$ or $\phat_{\bg}$ but does not exceed the asymptotic
minimax level by much. 

\begin{table}[htbp]
	\centering
	\begin{tabular}{|c|c|c|c|c|c|c|}
		\hline
		\textbf{Sparsity} & \textbf{r} & \textbf{Asymp-Theory} & \textbf{H-Plugin} & \textbf{C-Thresh} & \textbf{Bi-Grid} & \textbf{SS} \bigstrut\\
					\hline
					& 1     & 1.1513 & 120.4\% & 83.3\% & 86.6\% & 104.9\% \bigstrut\\
					\cline{2-7}          & 0.5   & 1.5351 & 173.6\% & 109.2\% & 105.5\% & 117.4\% \bigstrut\\
					\cline{2-7}    0.1   & 0.25  & 1.8421 & 278.5\% & 129.9\% & 125.5\% & 131.6\% \bigstrut\\
					\cline{2-7}          & 0.1   & 2.0933 & 588.1\% & 144.2\% & 156.0\% & 145.2\% \bigstrut\\
					\hline
					& 1     & 3.4539 & 109.1\% & 70.5\% & 70.8\% & 87.2\% \bigstrut\\
					\cline{2-7}          & 0.5   & 4.6052 & 162.1\% & 86.4\% & 84.6\% & 96.9\% \bigstrut\\
					\cline{2-7}    0.001 & 0.25  & 5.5262 & 267.6\% & 89.4\% & 95.9\% & 107.0\% \bigstrut\\
					\cline{2-7}          & 0.1   & 6.2798 & 582.8\% & 106.9\% & 114.0\% & 117.6\% \bigstrut\\
					\hline
					& 1     & 11.5129 & 123.9\% & 79.6\% & 78.9\% & 86.7\% \bigstrut\\
					\cline{2-7}    1E-10 & 0.5   & 15.3506 & 185.4\% & 88.0\% & 87.3\% & 94.2\% \bigstrut\\
					\cline{2-7}          & 0.25  & 18.4207 & 308.4\% & 94.6\% & 96.0\% & 100.3\% \bigstrut\\
					\cline{2-7}          & 0.1   & 20.9326 & 677.0\% & 101.8\% & 101.9\% & 106.0\% \bigstrut\\
					\hline
\end{tabular}%
\caption{\small \sl Numerical evaluation of the maximum risk
for the different univariate predictive densities over
$[-l,l]$ as the degree of sparsity $(\eta)$ and predictive
difficulty $r$ varies. Here, we have chosen $l=5\lambda$,
where $\lambda$ is defined in \eqref{eq:lambdadefs}. In
`Asymp-Theory' column we report the asymptotic minimax risk
$\lambda^2/(2r) $. In the other columns, we report the maximum risk of the estimators as quotients of the `Asymp-Theory' risk. }\label{numerical.risk.1}
\end{table}

\newpage

\newpage

\section{Future work}
 In this paper our results are based on known sparsity
 levels. Recently, computationally tractable Bayesian methods which
 adapts to unknown sparsity levels and possibly dense signals have
 been developed for point estimation
 \citep{carvalho2010horseshoe,bhattacharya2015dirichlet,rockova2014spike}. Using
 the approach of \citet{Johnstone04,Johnstone05} for point estimation,
 an interesting future direction will be to construct adaptive
 predictive density estimates.


\section{Appendix: Proof Details}\label{sec:append}

\subsection{Proof of the ``risk decomposition" Lemma~\ref{lem:decomp}}
Using \eqref{eq:dec}, write the Bayes predictive density as
\begin{equation} \label{eq:bpd}
p_D(y|x) 
= \frac{\int \phi(y|\mu,r) \phi(x-\mu) \pi(d \mu) }{\int
	\phi(x-\mu) \pi( d\mu)}
= \phi(y|0,r) \frac{N(x,y)}{D(x)},
\end{equation}
after rewriting numerator and denominator in the first ratio
respectively as
\begin{equation*}
\pi_0 \phi(y|0,r) \phi(x) N(x,y), \qquad \text{and} \qquad 
\pi_0 \phi(x) D(x).
\end{equation*}
After simple algebra, we find
\begin{equation}
N(x,y) = \int \exp \left\{ \mu \left(x + \frac{y}{r}\right) -
\frac{\mu^2}{2} \left(1+\frac{1}{r}\right) \right\} \frac{\pi(d \mu)}{\pi_0}  
\label{eq:Nxy}
\end{equation}
and $D(x)$ is analogous, but without terms in $y$ and $r$.
Note also that
\begin{equation*}
\E_\theta \log \bigg(\frac{\phi(Y|\theta,r)}{\phi(Y|0,r)}\bigg)
= \E_\theta \left[ \frac{\theta Y}{r} -
\frac{\theta^2}{2r} \right] 
= \frac{\theta^2}{2r}.
\end{equation*}
Hence, from (\ref{eq:bpd}) and the definition of predictive loss
\begin{equation*}
L(\theta,p_D(\cdot|x))
= \E_\theta \log \bigg(\frac{\phi(Y|\theta,r)}{p_D(Y|x)}\bigg)
= \frac{\theta^2}{2r} - \E_\theta \log N(x,Y) + \log D(x).
\end{equation*}

To obtain $\rho(\theta,\hat{p}_\pi)$, take expectation also over $X
\sim N(\theta,1)$. Since $Y \sim N(\theta,r)$ independently of $X$,
the random variable $X + Y/r \sim \mathcal{N}(\theta/v, 1/v)$ may be
expressed in the form $\theta/v + Z/\sqrt{v}$. 
Recalling the sparse prior form $\pi(d\mu) = (1 - \eta) \delta_0 +
\eta \nu$, we get
\begin{equation*}
N(X,Y) 
\stackrel{\mathcal{D}}{=} 1 + \frac{\eta}{1-\eta} \int \exp \left\{ \frac{\mu Z}{\sqrt v} 
+ \frac{\mu \theta}{v} -
\frac{\mu^2}{2v} \right \} \nu(d \mu)
= N_{\theta,v}(Z),
\end{equation*}
and the Lemma follows from the previous two displays.

\subsection{Proofs of the lemmas used in Section~\ref{sec:2-proofs}}\hspace{0.1cm}\\[1.5ex]
\textbf{Proof of Lemma \ref{lem:bounds}.} \ 
We do the easy lower bounds involving $N_{\theta,v}(Z)$ first. Indeed,
the bound for $\theta < \lambda$ follows 
just from $N_{\theta,v}(Z) \geq 1$. For $\theta = \lambda(\alpha_l + \omega) \geq \lambda$, by \eqref{eq:Al} we know:
\begin{align*}
 \E \log N_l &= \log c_1(\eta) + \hf \lambda^2 n(l,\omega), \text{ and, }  \\
 \E \log N_{l+1} &= \log c_1(\eta) + \hf \lambda^2 \{n(l+1,\omega)-2 \dot \alpha_l \, \alpha_{l+1}/v\}~.
\end{align*}
But $\log c_1(\eta) = \log c(\eta) - \log(1-\eta)^{-1} = O(1)$ as $\lambda \to \infty$. Hence, the proof of the lower bound is completed by using
\begin{equation*}
\E  \log N_{\theta,v}(Z)
\geq \max\{\E \log N_l, ~\E \log N_{l+1}\}.
\end{equation*}
\par
The proof of the upper bound on $\E\log D_{\theta}(Z)$ is more
involved,
and we first outline the approach. 
From \eqref{eq:Dtheta} and $1+x+y<(1+x)(1+y/x)$, we have
\begin{equation}\label{eq:lem2-temp1}
\log D_\theta(Z) 
\leq \log(1+ D_l) + \log(1+ \check D_l), 
\end{equation}
where we set 
$\check D_l  = \sum_{i \notin \{0,l\}} D_i/D_l$.
Henceforth in the proof, we make the choice $l = l(\theta)$ except
that when $0 \leq \theta < \mu_1$ we set $l = 1$.

For the first term (henceforth we call it the main term) in \eqref{eq:lem2-temp1} we will show
\begin{equation}
\label{eq:mainterm}
\E \log (1+D_l) \leq
\begin{cases}
\hf \lambda^2
d^+(l,\omega) + O(\lambda) & \quad \text{for } l \geq 1 \\
O(1)                       & \quad \text{if } 0 \leq \theta < \lambda

\end{cases}
\end{equation}
with $O(\lambda)$ being uniform in $l$. For the other term in
\eqref{eq:lem2-temp1} we will show that it is $O(\lambda)$ for all $l$
(and so, henceforth we call it the remainder term). For that purpose,
we write $D_{i,l} = D_i/D_l$ and decompose
\begin{equation*}
\check D_l 
= \sum_{k=1}^\infty D_{l+k,l}
+ \sum_{k=1}^{l-1} D_{l-k,l} 
+ \sum_{j=1}^\infty D_{-j,l}.
\end{equation*}
We have $ D_{-j}   \stackrel{\mathcal{D}}{=}  D_j \exp \{ - 2\mu_j \theta \}
\leq D_j $ since $\mu_j = - \mu_j, \pi_{-j} = \pi_j$ and $\mathcal{L}(Z)$ is
symmetric. Hence
\begin{equation*}
\sum_{j=1}^\infty D_{-j,l}
\stackrel{\mathcal{D}}{\leq} \sum_{j=1}^\infty D_{j,l}
= \sum_{k=1}^{l-1} D_{l-k,l} + 1 + \sum_{k=1}^\infty D_{l+k,l}.
\end{equation*}
Now use the elementary inequality $\log(1+ \sum \gamma_m )
\leq \sum \log (1+\gamma_m)$ to obtain
\begin{equation}
\label{eq:rem}
\E \log(1+ \check D_l)
\leq 2 \E  \log \Bigl( 1 + \sum_{k=1}^\infty D_{l+k,l}
\Bigr) + \log 2 + 
2 \E \log \Bigl( 1 + \sum_{k=1}^{l-1} D_{l-k,l}  \Bigr) .
\end{equation}
We will later show that the two main right side terms are each $O(\lambda)$.
This concludes the outline; we now turn to detailed analysis.

\bigskip
\textit{The Main term in \eqref{eq:lem2-temp1}.}
We first dispose of the case $0 \leq \theta <
\lambda$. 
From \eqref{eq:Aj} and \eqref{eq:piratio},
\begin{equation*}
D_1 = c_1(\eta) \exp \{ \lambda Z + \lambda \theta - \hf
\lambda^2(2+r^{-1}) \}.
\end{equation*}
Since $\theta < \lambda$ and $c_1(\eta) < (1-\eta)^{-1}$, and using
$\log(1+x) \leq \log 2 + (\log x)_+$, 
\begin{equation*}
\log(1 + D_1) \leq \log 2 + \log (1-\eta)^{-1} + \lambda Z
\end{equation*}
and hence $\E \log (1+D_1) \leq O(1)$. 

Now suppose that $\theta = \lambda(\alpha_l + \omega) \geq \lambda$
and use representation 
\eqref{eq:Al} for $D_l$. Abbreviating $\hf \lambda^2 d(l,\omega)$ 
as $d_{l \omega}$, we obtain
\begin{align*}
\E \log(1+D_l)
& = \E \log D_l + \E \log(1+D_l^{-1}) \notag \\
& = \log c(\eta) + \log(1-\eta)^{-1}  + d_{l \omega} + \log 2 + \E
(\log D_l^{-1})_+.
\intertext{Symmetry of  $\mathcal{L}(Z)$ about $0$ implies that
	$\log D_l^{-1} \stackrel{\mathcal{D}}{=} 
	- \log c(\eta) + \log(1-\eta) + \mu_l Z - d_{l \omega}$.
	By inspection the normalization constant $c(\eta) < 1$, so}
\E  (\log D_l^{-1})_+ 
& \leq - \log c(\eta) + \E (\mu_l Z - d_{l \omega})_+.
\end{align*}
From the previous two displays and  $\log (1-\eta)^{-1}
= O(\eta)$, we have
\begin{equation}
\label{eq:c}
\E \log(1+D_l)  = d_{l \omega} + \E (\mu_l Z - d_{l \omega})_+ + O(1).
\end{equation}
We now bound the expectation on the right side. 
Consider first those $l$ for which 
$\alpha_l \leq 2 + r^{-1}$ and thus
$\mu_l \leq (2 + r^{-1}) \lambda$. 
Noting that
\begin{equation*}
\E (\mu_l Z - d_{l \omega})_+
\leq - d_{l \omega} I \{ d_{l \omega} \leq 0 \} + \mu_l\, \E Z_+,
\end{equation*}
we then conclude that
\begin{equation*}
d_{l \omega} + \E (\mu_l Z - d_{l \omega})_+ 
\leq (d_{l \omega})_+ + (2+r^{-1}) \phi(0) \lambda.
\end{equation*}

Now consider the remaining $l$, with $\alpha_l \geq 2 + r^{-1}$,
for which we claim that
\begin{equation}
\label{eq:alphalbd}
\alpha_l^2 - \beta_l - r^{-1} \geq \hf \alpha_l^2.  
\end{equation}
We verify this via the equivalent form 
$ \alpha_l^2 - 2 \beta_l \geq 2 r^{-1}$. 
We have
\begin{equation*}
\frac{d}{dl} (\alpha_l^2 - 2 \beta_l)
= 2(\alpha_l \dot \alpha_l - \dot \beta_l) > 0
\end{equation*}
for both $1 \leq l \leq K$ and $l > K$.
So if $l_0 \in \RR_+$ satisfies $\alpha_{l_0} = 2 + r^{-1}$, we have
\begin{equation*}
\alpha_l^2 - 2 \beta_l
\geq \alpha_{l_0}^2 - 2 \beta_{l_0} 
\geq \alpha_{l_0}^2 - 2 \alpha_{l_0} 
= (2 + r^{-1})r^{-1} \geq 2 r^{-1}.
\end{equation*}

Since $\omega \geq 0$, we have from \eqref{eq:n} and \eqref{eq:alphalbd},
\begin{equation*}
d_{l \omega} \geq \hf \lambda^2[\alpha_l^2 - \beta_l - r^{-1}]
\geq \tfrac{1}{4} (\lambda \alpha_l)^2 
=    \tfrac{1}{4} \mu_l^2.
\end{equation*}
From the bound
$\E (Z-x)_+ \leq \phi(x)/x^2$ we calculate
\begin{equation*}
\E (\mu_l Z - d_{l \omega})_+
\leq \mu_l \E(Z -  \mu_l/4)_+
\leq 16 \frac{\phi(\mu_l/4)}{ \mu_l} \leq M_b,
\end{equation*}
uniformly in $\lambda \geq 1$ and $l$ such that $\alpha_l \geq 2 +
r^{-1}$. 
Combining the two cases with \eqref{eq:c}, we have proven the
bound \eqref{eq:mainterm} on the first term of
\eqref{eq:lem2-temp1}.

\bigskip
We turn now to bounding the remainder \eqref{eq:rem}.
This depends on the decay between successive terms $D_j$, so we start
by using \eqref{eq:altDrep} to derive a useful representation for $D_{j+1}/D_j$.
Indeed, using $\mu_j = \lambda \alpha_j$ and $\theta = \lambda(\alpha_l +
\omega)$, we define 
\begin{align*}
\Delta_j
= \Delta(j; l, \omega)
& = (2/\lambda^2) [G(\mu_{j+1};\theta) - G(\mu_j;\theta)] \\
& = \dot \alpha_j [ \alpha_{j+1} + \alpha_j - 2 \alpha_l - 2 \omega]
+ \dot \beta_j
\end{align*}
and arrive at, for $j \geq 1$,
\begin{equation}
\label{eq:rep}
\frac{D_{j+1}}{D_j} 
= \exp \{ \lambda \dot \alpha_j Z - \hf \lambda^2
\Delta_j \}.
\end{equation}

We now show that $l \to \Delta_l$ achieves a maximum at $l = 0$.  This
will also verify the claim in Section 2 that $j \to G(\mu_j;\theta)$
is minimized at $j = l(\theta)$ for each
$\theta \in [\mu_l,\mu_{l+1})$. The argument splits into two largely
parallel cases. 

Suppose first that $j \geq l$, so that $j = l+k$ for $k \geq 0$.
Using $\alpha_l + \omega \leq \alpha_{l+1}$, then
$\dot \beta_{l+k} = \dot \alpha_{l+k}^2$ and finally
$\dot \alpha_{l+k} + \alpha_{l+k} = \alpha_{l+k+1} $,
we have
\begin{equation}
\label{eq:kpos}
\begin{split}
\Delta_{l+k} 
&\geq \dot \alpha_{l+k} [\alpha_{l+k+1} + \alpha_{l+k} - 2
\alpha_{l+1}] + \dot \alpha_{l+k}^2 \\
& = 2 \dot \alpha_{l+k} (\alpha_{l+k+1} - \alpha_{l+1})
\geq 0,
\end{split}
\end{equation}
with strict inequality when $k \geq 1$.

Suppose now that $j < l$, so that $j = l-k-1$ for $k \geq 0$.
Using $\alpha_l + \omega \geq \alpha_{l}$, then
$\dot \beta_{l-k-1} = \dot \alpha_{l-k-1}^2$ and finally
$\dot \alpha_{l-k-1} + \alpha_{l-k-1} = \alpha_{l-k} $,
we have
\begin{align*}
\Delta_{l-k-1} 
&\leq \dot \alpha_{l-k-1} [\alpha_{l-k} + \alpha_{l-k-1} - 2
\alpha_{l}] + \dot \alpha_{l-k-1}^2 \\
& = 2 \dot \alpha_{l-k-1} (\alpha_{l-k} - \alpha_{l})
\leq 0,
\end{align*}
with strict inequality when $k \geq 1$. 

\bigskip
As final preparation, we record a useful bound whose proof is provided at the end this proof.

\begin{lemma}\label{pf.D2.lem.2}
	If $a_1, a_2, \ldots$ are positive, then for each $n \geq 1$,
	\begin{equation}
	\label{eq:pf.D2.lem.2-2}
	\log \bigg(1+\sum_{k=1}^{n+1} a_k \bigg)
	< \log(1+ a_1) + \sum_{k=1}^n \frac{a_{k+1}}{a_k}. 
	\end{equation}
\end{lemma}


We next concentrate on bounding the \textit{first term of \eqref{eq:rem}.} \ 
Use \eqref{eq:pf.D2.lem.2-2} with $a_k = D_{l+k}/D_l$ and
$\log(1+a_1) \leq \log 2 + (\log a_1)_+$ to write
\begin{equation}
\label{eq:first}
\E \log \bigg(1+ \sum_{k=1}^\infty D_{l+k,l} \bigg)
\leq \log 2 + \E \Bigl(\log \frac{D_{l+1}}{D_l} \Bigr)_+ 
+ \E \bigg\{\sum_{k=1}^\infty  \frac{D_{l+k+1}}{D_{l+k}}\bigg\}. 
\end{equation}
In \eqref{eq:rep} with $j=l$, we have seen that $\Delta_l \geq 0$ and so
\begin{equation*}
\E \Bigl(\log \frac{D_{l+1}}{D_l} \Bigr)_+  \leq \lambda \dot
\alpha_l \,\E Z_+ \leq \lambda \phi(0).
\end{equation*}
When $j=l+k$, observe from 
\eqref{eq:kpos} that
$\Delta_{l+k} \geq 2 \dot \alpha_{l+k}^2 + 2 \dot
\alpha_{l+k}(\alpha_{l+k} - \alpha_{l+1})$.
From \eqref{eq:rep}, now with $j=l+k$ for $k \geq 1$,
\begin{align*}
\E \bigg\{\frac{D_{l+k+1}}{D_{l+k}}\bigg\}
& = \exp \{ \hf \lambda^2 [\dot \alpha_{l+k}^2 - \Delta_{l+k}] \} \\
& \leq \exp \{ - \hf \lambda^2 [\dot \alpha_{l+k}^2 + 2 \dot
\alpha_{l+k}(\alpha_{l+k} - \alpha_{l+1})] \} \\
& \leq \exp \{ - \hf \lambda^2 b^2 - \lambda^2 b^2 (k-1) \},
\end{align*}
so that the right side of \eqref{eq:first} is 
$O(\lambda) + O(e^{-\lambda^2 b^2/2}) = O(\lambda)$. 

\bigskip
\textit{Second term of \eqref{eq:rem}.} \ 
Now use \eqref{eq:pf.D2.lem.2-2} with $a_k = D_{l-k}/D_l$:
\begin{equation}
\label{eq:second}
\E \log \bigg(1+ \sum_{k=1}^{l-1} D_{l-k,l} \bigg)
\leq \log 2 + \E \Bigl(\log \frac{D_{l-1}}{D_l} \Bigr)_+ 
+ \E \bigg\{\sum_{k=1}^{l-2}  \frac{D_{l-k-1}}{D_{l-k}}\bigg\}. 
\end{equation}
In \eqref{eq:rep} with $j=l-1$, we have seen that $\Delta_{l-1} \leq 0$ and so
\begin{equation*}
\E \Bigl(\log \frac{D_{l-1}}{D_l} \Bigr)_+  \leq \lambda \dot
\alpha_{l-1} \E Z_+ \leq \lambda \phi(0).
\end{equation*}
From \eqref{eq:rep}, now with $j=l-k-1$,
\begin{align*}
\E \bigg\{\frac{D_{l-k-1}}{D_{l-k}}\bigg\}
& = \E \{\exp \{ - \lambda \dot \alpha_{l-k-1} Z + \hf \lambda^2
\Delta_{l-k-1} \}\} \\
& \leq \exp \{ \hf \lambda^2 \dot \alpha_{l-k-1} [\dot \alpha_{l-k-1}
+ 2(\alpha_{l-k} - \alpha_l)] \},
\end{align*}
as $\Delta_{l-k-1} \leq 2 \dot \alpha_{l-k-1} (\alpha_{l-k} - \alpha_{l})$. 
Since  $j \to \dot \alpha_j$ is increasing and $\dot \alpha_j \geq b$,
\begin{align*}
\dot \alpha_{l-k-1} + 2(\alpha_{l-k}-\alpha_l)
& \leq \dot \alpha_{l-k}  - 2(\alpha_l -
\alpha_{l-k+1})  -  2\dot \alpha_{l-k} \\
& \leq -b -2(k-1)b.
\end{align*}
Using  $\dot \alpha_{l-k-1} \geq b$ again, we conclude that
\begin{equation*}
\E \bigg\{\sum_{k=1}^{l-2}  \frac{D_{l-k-1}}{D_{l-k}} \bigg\}
\leq \sum_{k=1}^\infty \exp \{ - \hf \lambda^2 b^2 - \lambda^2 b^2
(k-1) \} 
= O(e^{-\lambda^2 b^2/2}).
\end{equation*}
Thus, we have proved the desired bound on the second term.  This completes the proof of the lemma.\medskip
\\
\noindent \textbf{Proof of Lemma~\ref{pf.D2.lem.2}.}
	We use induction. The bounds 
	$1+x+y<(1+x)(1+y/x)$ and $\log(1+x)<x$, valid for positive $x,y$,
	establish the case $n=1$. For general $n$, 
	let $a_k' = a_{k+1}/(1+a_1)$ for $k = 1, \ldots, n$. Then
	\begin{align*}
	\log \bigg(1+ \sum_{k=1}^{n+1}a_k \bigg)
	& = \log(1+a_1) + \log \bigg(1 + \sum_{k=1}^n a_k' \bigg) \\
	& < \log(1+a_1) + \log(1 + a_1') + \sum_{k=1}^{n-1}
	\frac{a_{k+1}'}{a_k'}  \\
	& < \log(1+a_1) + \frac{a_2}{a_1} + \sum_{k=2}^{n}
	\frac{a_{k+1}}{a_{k}},
	\end{align*}
	where the inequalities use the cases $n-1$ and $n=1$ in turn.\\[1ex]

\noindent \textbf{Proof of Lemma~\ref{lem.pf.D3}.}
%
We modify some of the methods used in Lemma~\ref{lem:bounds} to
incorporate now $v = r(1+r)^{-1} \in (0,1)$. 
With $N_j$ defined as in \eqref{eq:N} and \eqref{eq:Nj}, and arguing
as around \eqref{eq:lem2-temp1}, we bound
\begin{equation*}
  \log N_{\theta,v}(Z) \leq \log (1+N_1) + \log (1+ \check N_1),
\end{equation*}
with $\check N_1=\sum_{j \notin \{0,1\}} N_{j,1}$ and 
$N_{j,1}=N_j/N_1$.
The desired control of the main term is easier here than in
Lemma~\ref{lem:bounds}: for $l=1$, 
\begin{align*}
\ex \log (1+N_1) & \leq \log 2 + \ex (\log N_1)_+ \\
&\leq \log 2 + \log (1-\eta)^{-1} + v^{-1/2} \lambda \, \ex|Z| + 2^{-1}\lambda^2 n(1,\omega)\\
&\leq 2^{-1}\lambda^2 n(1,\omega) + O(\lambda).
\end{align*}

Turn now to the remainder term. 
Since $N_{-j}\stackrel{D}{=} N_j \exp(-2 v^{-1/2}\mu_j\theta)$, 
we may argue as before to obtain the analog of \eqref{eq:rem}:
\begin{equation*}
  \ex\{\log(1+\check N_1)\} \leq \log 2 + 2 \, \mathrm{Rem}(\lambda),
\end{equation*}
where, using Lemma~\ref{pf.D2.lem.2} with 
$a_k = N_{k+1,1}$ and setting $\check{R}_k = N_{k+1}/N_k$,
\begin{equation*}
  \mathrm{Rem}(\lambda) 
    = \ex \log \Big(1 + \sum_{k=1}^\infty N_{k+1,1} \Big)
    \leq \log 2 + \ex (\log \check{R}_1)_+ 
                + \sum_{k=1}^\infty \ex \, \check{R}_{k+1}.
\end{equation*}
Using definition \eqref{eq:Nj} and then \eqref{eq:piratio} and
\eqref{eq:def.spacing}, we obtain
\begin{align*}
  \check{R}_k
   & = (\pi_{k+1}/\pi_k) \exp
     \{v^{-1/2}(\mu_{k+1}-\mu_k)(Z+\theta) - \hf
     v^{-1}(\mu_{k+1}^2-\mu_k^2)\} \\
   & = \exp \{ \lambda v^{-1/2} \dot \alpha_k(Z+\theta) 
              -\hf \lambda^2 [ v^{-1}\dot
     \alpha_k(\alpha_{k+1}+\alpha_k) - \dot \beta_k] \} \\ 
   & = \exp \{ \lambda v^{-1/2} \dot \alpha_k Z - \hf \lambda^2 v^{-1}
     \gamma_k(\omega,v) \},     
\end{align*}
where we put $\theta = \lambda(1+\omega)$ and used
$\dot \beta_k = \dot \alpha_k^2$ and $\alpha_{k+1} = \dot \alpha_k + 
\alpha_k$ to write
\begin{align*}
  \gamma_k(\omega,v)
  & = -2v^{1/2} \dot \alpha_k(1+\omega) + \dot \alpha_k^2 + 2 \dot
    \alpha_k \alpha_k + v \dot \alpha_k^2 \\
  & = \dot \alpha_k[(1+v)\dot \alpha_k + 2(\alpha_k-\sqrt{v}(1+\omega))].
\end{align*}

Since $b \leq 1$, we necessarily have $K - 1 = \lceil 2b^{-3/2} \rceil
\geq 2$, and therefore $\alpha_1 = 1, \alpha_2 = 1+b$ and $\dot \alpha_1 =
b$. 
Consequently, at $k = 1$ we have
\begin{equation*}
  \gamma_1(\omega,v) 
    = b[b + bv + 2 - 2 \sqrt{v}(1+\omega)]
    \geq b(bu^2+2u),
\end{equation*}
since $\omega\leq b$ and we have set $u = 1 - \sqrt v \geq 0$.
We arrive at
\begin{equation*}
  \log \check{R}_1 
    \leq \lambda v^{-1/2} b Z - \hf \lambda^2 v^{-1} bu(bu+2),
\end{equation*}
and therefore $\ex (\log \check{R}_1)_+ \leq \lambda v^{-1/2} b \ex
Z_+ = O(\lambda).$

For $k\geq2$, we write 
$\log \ex \check{R}_k = - \hf \lambda^2 v^{-1}(\gamma_k-\dot
\alpha_k^2)$, and from $\dot \alpha_k \geq b$,
\begin{equation*}
  \gamma_k-\dot \alpha_k^2
    = \dot \alpha_k[v\dot \alpha_k + 2(\alpha_k-\sqrt{v}(1+\omega))]
   \geq b^2(v+2(k-2)).
\end{equation*}
since $\alpha_k - \sqrt{v}(1+\omega) \geq \alpha_k - \alpha_2 \geq (k-2)b$.
This entails
\begin{equation*}
  \sum_{k=2}^\infty \ex \check{R}_{k} 
    \leq e^{-b^2 \lambda^2/2} 
         \sum_{k=0}^\infty \exp (-v^{-1}\lambda^2 b^2 k) 
    = O(\lambda).
\end{equation*}
The last two paragraphs show that $\textrm{Rem}(\lambda) = O(\lambda)$
and complete the proof.

\subsection{Detailed proof of the lemmas and the inequalities used Section~\ref{sec-2-3}}\hspace{1mm}\\
	
\noindent \textbf{Proof of Lemma 3.1.}
	From \eqref{eq:mean}, we find
	\begin{equation*}
	\E N_{\theta,v}^S(Z)
	\leq 1+\frac{\eta}{1-\eta} \frac{1}{2l} \int_{-\infty}^l \exp \Big(
	\frac{\mu \theta}{v} \Big) \nu(d \mu)
	\leq 1 + \frac{\eta}{(1-\eta)} \frac{e^w}{2w},
	\end{equation*}
	for $w = \theta l/v \geq 1$. 
	From this and Jensen's inequality \eqref{eq:jensen}, we obtain
	\begin{equation*}
	\E \log N_{\theta,v}^S(Z) \leq w + \log \left\{ e^{-w} + \frac{\eta}{(1-\eta)}
	\frac{1}{2w} \right\} \leq w, 
	\end{equation*}
	since the term in braces is bounded by $1/e + 1/2 \leq 1$. \\[1ex]

\noindent\textbf{Proof of (\ref{eq:Phi-bd}).} Use $v \leq 1$, then $0 \leq \theta \leq l$ and finally $l \geq 1$ to
conclude 
\begin{equation*}
\Phi_{l,v} \geq \Phi_{l,1} 
\geq \Phi(-Z) - \Phi(-l-Z) \geq \Phi(-Z) - \Phi(-1-Z).
\end{equation*}
Now use symmetry of $Z$ and then Jensen's inequality to get
\begin{align*}
\E \log \Phi_{l,v}
\geq \E \log \int_{Z}^{Z+1} \phi(s)\, ds 
&\geq \E \int_{Z}^{Z+1} \log \phi(s)\, ds\\
& = \log \phi(0) - 2/3.
\end{align*}

\bigskip
\noindent \textbf{Proof of (\ref{eq:D-bd}).} From the definition (\ref{eq:Nrep}) with $v = 1$ and using
$\Phi_{l,1} \leq 1$, 
$$\log D_{\theta}^S(Z) \leq  \log 2 + [\log c(\eta) + (\theta +
Z)^2/2]_+ .$$ 
From the upper bound in (\ref{eq:ceta-bds}), we have
$$\E\log D_{\theta}^S(Z) \leq  2^{-1}\E [(\theta+Z)^2-v^{-1}
\lambda^2]_+ + O(1).$$ 

Now if  $\theta <\tilde{\lambda}$, then
\begin{align*}
\E [(\theta+Z)^2-v^{-1} \lambda^2]_+ &\leq [\theta^2-v^{-1}
\lambda^2]_+ + 2 |\theta| \,
\E Z_+ + \E Z^2 \\
& \leq [\tilde{\lambda}^2-v^{-1}\lambda^2]_+ + 2  \, \tilde{\lambda}
\, \phi(0) + 1 
\end{align*}
which suffices for the first bound.

If $\theta \geq \tilde{\lambda}$,  we put
$W =(\theta+Z)^2$ and $c = \lambda^2/v$
and apply the inequality
$\E (W-c)_+ \leq \E(W-c) + c P(W<c)$, valid for $W \geq 0$.
For all $\theta \geq \tilde \lambda$, we have
$\{ |\theta + Z | \leq \lambda/\sqrt{v} \} \subset \{ Z < - \sqrt{2
	\log \lambda} \}$ and so
\begin{align*}
\E [(\theta+Z)^2-v^{-1} \lambda^2]_+ 
&\leq \E (\theta+Z)^2-v^{-1} \lambda^2 
+ v^{-1} \lambda^2P( Z > \sqrt{ 2 \log \lambda})\\
& \leq \theta^2 -v^{-1} \lambda^2 + O(\lambda) \text{ as } \lambda \to \infty~.
\end{align*}
where the last inequality uses the Mills ratio bound $ P( Z > x) \leq
x^{-1} \phi(x)$. 



\section*{Acknowledgements}
\label{sec:acknowledgements}
GM was supported in part by the Zumberge individual award from the University of Southern California's James H. Zumberge faculty research and innovation fund. IMJ was supported in part by NSF DMS 1407813 and 1418362 and  thanks the Australian National
University for hospitality while working on this paper.


\bibliographystyle{natbib}
\bibliography{pred-inf}

\def\cprime{$'$}
\begin{thebibliography}{}

\bibitem[Aitchison and Dunsmore(1975)Aitchison and Dunsmore]{Aitchison-book}
Aitchison, J. and Dunsmore, I.~R. (1975).
\newblock {\em Statistical prediction analysis\/}.
\newblock Cambridge University Press, Cambridge.

\bibitem[Aslan(2006)Aslan]{Aslan06}
Aslan, M. (2006).
\newblock Asymptotically minimax {B}ayes predictive densities.
\newblock {\em Ann. Statist.}, {\bf 34}(6), 2921--2938.

\bibitem[Bhattacharya {\em et~al.}(2015)Bhattacharya, Pati, Pillai, and
  Dunson]{bhattacharya2015dirichlet}
Bhattacharya, A., Pati, D., Pillai, N.~S., and Dunson, D.~B. (2015).
\newblock Dirichlet--laplace priors for optimal shrinkage.
\newblock {\em Journal of the American Statistical Association\/}, {\bf
  110}(512), 1479--1490.

\bibitem[Brown {\em et~al.}(2008)Brown, George, and Xu]{Brown08}
Brown, L.~D., George, E.~I., and Xu, X. (2008).
\newblock Admissible predictive density estimation.
\newblock {\em Ann. Statist.}, {\bf 36}(3), 1156--1170.

\bibitem[Carvalho {\em et~al.}(2010)Carvalho, Polson, and
  Scott]{carvalho2010horseshoe}
Carvalho, C.~M., Polson, N.~G., and Scott, J.~G. (2010).
\newblock The horseshoe estimator for sparse signals.
\newblock {\em Biometrika\/}, pages 465--480.

\bibitem[Fourdrinier {\em et~al.}(2011)Fourdrinier, Marchand, Righi, and
  Strawderman]{Fourdrinier11}
Fourdrinier, D., Marchand, {\'E}., Righi, A., and Strawderman, W.~E. (2011).
\newblock On improved predictive density estimation with parametric
  constraints.
\newblock {\em Electron. J. Stat.}, {\bf 5}, 172--191.

\bibitem[Geisser(1993)Geisser]{Geisser-book}
Geisser, S. (1993).
\newblock {\em Predictive inference\/}, volume~55 of {\em Monographs on
  Statistics and Applied Probability\/}.
\newblock Chapman and Hall, New York.
\newblock An introduction.

\bibitem[George and McCulloch(1997)George and McCulloch]{george1997approaches}
George, E.~I. and McCulloch, R.~E. (1997).
\newblock Approaches for bayesian variable selection.
\newblock {\em Statistica sinica\/}, pages 339--373.

\bibitem[George and Xu(2008)George and Xu]{George08}
George, E.~I. and Xu, X. (2008).
\newblock Predictive density estimation for multiple regression.
\newblock {\em Econometric Theory\/}, {\bf 24}(2), 528--544.

\bibitem[George {\em et~al.}(2006)George, Liang, and Xu]{George06}
George, E.~I., Liang, F., and Xu, X. (2006).
\newblock Improved minimax predictive densities under {K}ullback-{L}eibler
  loss.
\newblock {\em Ann. Statist.}, {\bf 34}(1), 78--91.

\bibitem[George {\em et~al.}(2012)George, Liang, and Xu]{George12}
George, E.~I., Liang, F., and Xu, X. (2012).
\newblock From minimax shrinkage estimation to minimax shrinkage prediction.
\newblock {\em Statist. Sci.}, {\bf 27}(1), 82--94.

\bibitem[Ghosh {\em et~al.}(2008)Ghosh, Mergel, and Datta]{Ghosh08}
Ghosh, M., Mergel, V., and Datta, G.~S. (2008).
\newblock Estimation, prediction and the {S}tein phenomenon under divergence
  loss.
\newblock {\em J. Multivariate Anal.}, {\bf 99}(9), 1941--1961.

\bibitem[Hartigan(1998)Hartigan]{Hartigan98}
Hartigan, J.~A. (1998).
\newblock The maximum likelihood prior.
\newblock {\em Ann. Statist.}, {\bf 26}(6), 2083--2103.

\bibitem[Ishwaran and Rao(2005a)Ishwaran and Rao]{ishwaran2005spike-a}
Ishwaran, H. and Rao, J.~S. (2005a).
\newblock Spike and slab gene selection for multigroup microarray data.
\newblock {\em Journal of the American Statistical Association\/}, {\bf
  100}(471), 764--780.

\bibitem[Ishwaran and Rao(2005b)Ishwaran and Rao]{ishwaran2005spike}
Ishwaran, H. and Rao, J.~S. (2005b).
\newblock Spike and slab variable selection: frequentist and bayesian
  strategies.
\newblock {\em Annals of Statistics\/}, pages 730--773.

\bibitem[Johnstone(1994)Johnstone]{Johnstone94a}
Johnstone, I.~M. (1994).
\newblock On minimax estimation of a sparse normal mean vector.
\newblock {\em Ann. Statist.}, {\bf 22}(1), 271--289.

\bibitem[Johnstone(2013)Johnstone]{Johnstone-book}
Johnstone, I.~M. (2013).
\newblock Gaussian estimation: Sequence and wavelet models.
\newblock Version: 11 June, 2013. Available at
  \url{"http://www-stat.stanford.edu/~imj"}.

\bibitem[Johnstone and Silverman(2004)Johnstone and Silverman]{Johnstone04}
Johnstone, I.~M. and Silverman, B.~W. (2004).
\newblock Needles and straw in haystacks: empirical {B}ayes estimates of
  possibly sparse sequences.
\newblock {\em Ann. Statist.}, {\bf 32}(4), 1594--1649.

\bibitem[Johnstone and Silverman(2005)Johnstone and Silverman]{Johnstone05}
Johnstone, I.~M. and Silverman, B.~W. (2005).
\newblock Empirical {B}ayes selection of wavelet thresholds.
\newblock {\em Ann. Statist.}, {\bf 33}(4), 1700--1752.

\bibitem[Kobayashi and Komaki(2008)Kobayashi and Komaki]{Kobayashi08}
Kobayashi, K. and Komaki, F. (2008).
\newblock Bayesian shrinkage prediction for the regression problem.
\newblock {\em J. Multivariate Anal.}, {\bf 99}(9), 1888--1905.

\bibitem[Komaki(1996)Komaki]{Komaki96}
Komaki, F. (1996).
\newblock On asymptotic properties of predictive distributions.
\newblock {\em Biometrika\/}, {\bf 83}(2), 299--313.

\bibitem[Komaki(2001)Komaki]{Komaki01}
Komaki, F. (2001).
\newblock A shrinkage predictive distribution for multivariate normal
  observables.
\newblock {\em Biometrika\/}, {\bf 88}(3), 859--864.

\bibitem[Komaki(2006)Komaki]{komaki2006shrinkage}
Komaki, F. (2006).
\newblock Shrinkage priors for bayesian prediction.
\newblock {\em the Annals of Statistics\/}, pages 808--819.

\bibitem[Kubokawa {\em et~al.}(2013)Kubokawa, Marchand, Strawderman, and
  Turcotte]{kubokawa2013minimaxity}
Kubokawa, T., Marchand, {\'E}., Strawderman, W.~E., and Turcotte, J.-P. (2013).
\newblock Minimaxity in predictive density estimation with parametric
  constraints.
\newblock {\em Journal of Multivariate Analysis\/}, {\bf 116}, 382--397.

\bibitem[Kubokawa {\em et~al.}(2015)Kubokawa, Marchand, and
  Strawderman]{kubokawa2015predictive}
Kubokawa, T., Marchand, {\'E}., and Strawderman, W.~E. (2015).
\newblock On predictive density estimation for location families under
  integrated squared error loss.
\newblock {\em Journal of Multivariate Analysis\/}, {\bf 142}, 57--74.

\bibitem[Kubokawa {\em et~al.}(2017)Kubokawa, Marchand, Strawderman, {\em
  et~al.}]{kubokawa2017predictive}
Kubokawa, T., Marchand, {\'E}., Strawderman, W.~E., {\em et~al.} (2017).
\newblock On predictive density estimation for location families under
  integrated absolute error loss.
\newblock {\em Bernoulli\/}, {\bf 23}(4B), 3197--3212.

\bibitem[Mallows(1978)Mallows]{mall78}
Mallows, C. (1978).
\newblock Minimizing an integral.
\newblock {\em SIAM Review\/}, {\bf 20}(1), 183--183.

\bibitem[Maruyama and Ohnishi(2016)Maruyama and Ohnishi]{maruyama2016harmonic}
Maruyama, Y. and Ohnishi, T. (2016).
\newblock Harmonic bayesian prediction under alpha-divergence.
\newblock {\em arXiv preprint arXiv:1605.05899\/}.

\bibitem[Matsuda and Komaki(2015)Matsuda and Komaki]{matsuda2015singular}
Matsuda, T. and Komaki, F. (2015).
\newblock Singular value shrinkage priors for bayesian prediction.
\newblock {\em Biometrika\/}, {\bf 102}(4), 843--854.

\bibitem[Mitchell and Beauchamp(1988)Mitchell and
  Beauchamp]{mitchell1988bayesian}
Mitchell, T.~J. and Beauchamp, J.~J. (1988).
\newblock Bayesian variable selection in linear regression.
\newblock {\em Journal of the American Statistical Association\/}, {\bf
  83}(404), 1023--1032.

\bibitem[Mukherjee(2013)Mukherjee]{Mukherjee-thesis}
Mukherjee, G. (2013).
\newblock {\em Sparsity and Shrinkage in Predictive Density Estimation\/}.
\newblock Ph.D. thesis, Stanford University.

\bibitem[Mukherjee and Johnstone(2015)Mukherjee and Johnstone]{Mukherjee-15}
Mukherjee, G. and Johnstone, I.~M. (2015).
\newblock Exact minimax estimation of the predictive density in sparse gaussian
  models.
\newblock {\em Annals of Statistics\/}, {\bf 43}(3), 937.

\bibitem[O'Hara {\em et~al.}(2009)O'Hara, Sillanp{\"a}{\"a}, {\em
  et~al.}]{o2009review}
O'Hara, R.~B., Sillanp{\"a}{\"a}, M.~J., {\em et~al.} (2009).
\newblock A review of bayesian variable selection methods: what, how and which.
\newblock {\em Bayesian analysis\/}, {\bf 4}(1), 85--117.

\bibitem[Park and Casella(2008)Park and Casella]{park2008bayesian}
Park, T. and Casella, G. (2008).
\newblock The bayesian lasso.
\newblock {\em Journal of the American Statistical Association\/}, {\bf
  103}(482), 681--686.

\bibitem[Rockov{\'a}(2017)Rockov{\'a}]{rockova2015bayesian}
Rockov{\'a}, V. (2017).
\newblock Bayesian estimation of sparse signals with a continuous
  spike-and-slab prior.
\newblock {\em Annals of Statistics\/}.

\bibitem[Ro{\v{c}}kov{\'a} and George(2014)Ro{\v{c}}kov{\'a} and
  George]{rovckova2014negotiating}
Ro{\v{c}}kov{\'a}, V. and George, E.~I. (2014).
\newblock Negotiating multicollinearity with spike-and-slab priors.
\newblock {\em Metron\/}, {\bf 72}(2), 217--229.

\bibitem[Ro{\v{c}}kov{\'a} and George(2016)Ro{\v{c}}kov{\'a} and
  George]{rockova2014spike}
Ro{\v{c}}kov{\'a}, V. and George, E.~I. (2016).
\newblock The spike-and-slab lasso.
\newblock {\em Journal of the American Statistical Association\/}.

\bibitem[Xu and Liang(2010)Xu and Liang]{Xu10}
Xu, X. and Liang, F. (2010).
\newblock Asymptotic minimax risk of predictive density estimation for
  non-parametric regression.
\newblock {\em Bernoulli\/}, {\bf 16}(2), 543--560.

\bibitem[Xu and Zhou(2011)Xu and Zhou]{Xu12}
Xu, X. and Zhou, D. (2011).
\newblock Empirical bayes predictive densities for high-dimensional normal
  models.
\newblock {\em J. Multivariate Analysis\/}, {\bf 102}(10), 1417--1428.

\bibitem[Yano and Komaki(2017)Yano and Komaki]{yano2017information}
Yano, K. and Komaki, F. (2017).
\newblock Information criteria for prediction when the distributions of current
  and future observations differ.
\newblock {\em Statistica Sinica\/}, {\bf 27}, 1205--1223.

\end{thebibliography}
\newpage
\section{Supplementary Risk plots}
\begin{figure}[h!] 
	\includegraphics[width=\textwidth]{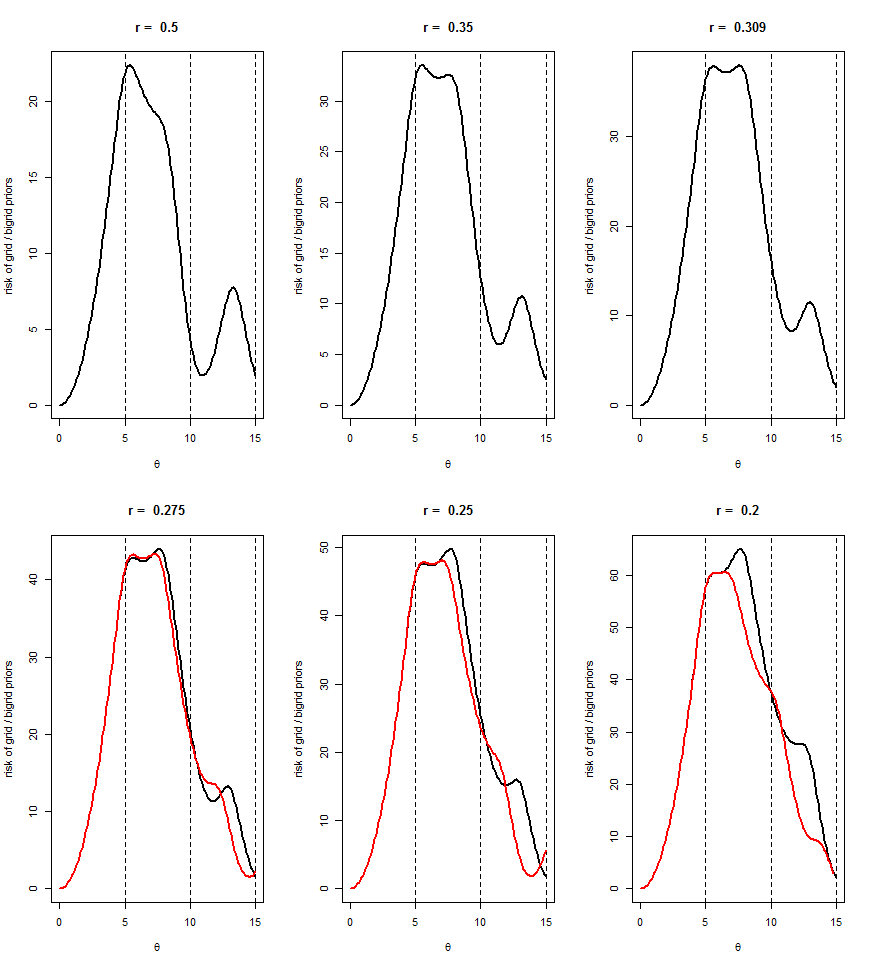}
	\caption{\small \sl Plot of the risk of $\phat_{\g}$ (in black) and $\phat_{\bg}$ (in red for $r< r_0=(\sqrt{5}-1)/4\approx 0.309$) as the parameter $\theta$ varies over $\mathbb{R}^+$. Here, $\lambda=5$ and $\eta=\exp(-\lambda^2/(2v))$.  The vertical line denotes multiples of $\lambda$ along the x-axis.} \label{fig-1a}
\end{figure}

\begin{figure}[h!] 
	\includegraphics[width=0.9\textwidth]{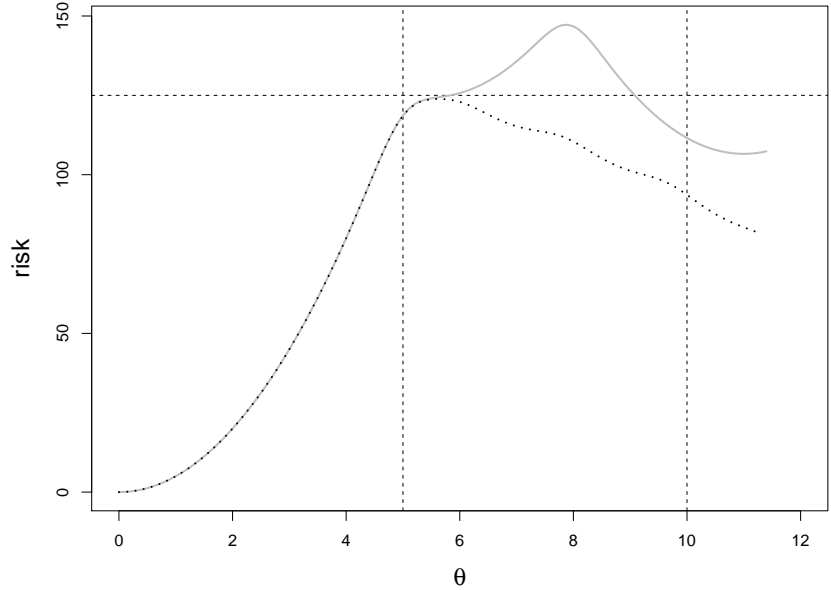}
	\caption{\small \sl In black dotted line we have the plot of $\phat_{\bg}$ and in gray we have plot of $\phat_{\g}$. The vertical lines mark $\lambda$ and $2\lambda$ while the horizontal line represents the bench mark $\lambda^2/(2 r)$. Here, $\lambda=5$, $r=0.1$ and $\eta=\exp(-\lambda^2/2)$. The maximal risk of $\phat_{\g}$ is 1.18 times the bench mark value. The risk of $\phat_{\bg}$ is controlled below the bench mark value.}\label{fig-6}
\end{figure}

\hspace{1cm}\\[1ex]
Here, through numerical evaluations we provide further insights on the risk properties for the  \textsf{pdes} considered in the previous sections. Figure~\ref{fig-1} below shows the risk of the grid prior is well controlled below the desired limit when $r=1$. Also, the plot reveals that the risk function exhibits periodicity for sufficiently large parametric values (in this case as $|\theta| \geq 2 \lambda$) with a period of $\lambda$. In Figure~\ref{fig-1a} we have the risk plots of the grid and bi-grid priors for different values of $r$. The bi-grid comes in play for $r< r_0\approx 0.309$. From the plots we see that for $|\theta| \in [\lambda, 2\lambda]$, the risk function for the grid prior is roughly decreasing in $|\theta|$ for large value of $r$ but is bi-modal as $r$ decreases towards $r_0$. At $r_0$ its two peaks are equal in height and as $r$ decreases further the gap between the maximal risk of the grid prior and the bi-grid prior widens. In Figure~\ref{fig-6}, we exhibit a scenario where the \textsf{pde} based on the grid prior is no longer optimal and its risk is  far dominated by the bi-grid prior based \textsf{pde}.
   
\begin{figure}[h!] 
	\includegraphics[width=0.9\textwidth]{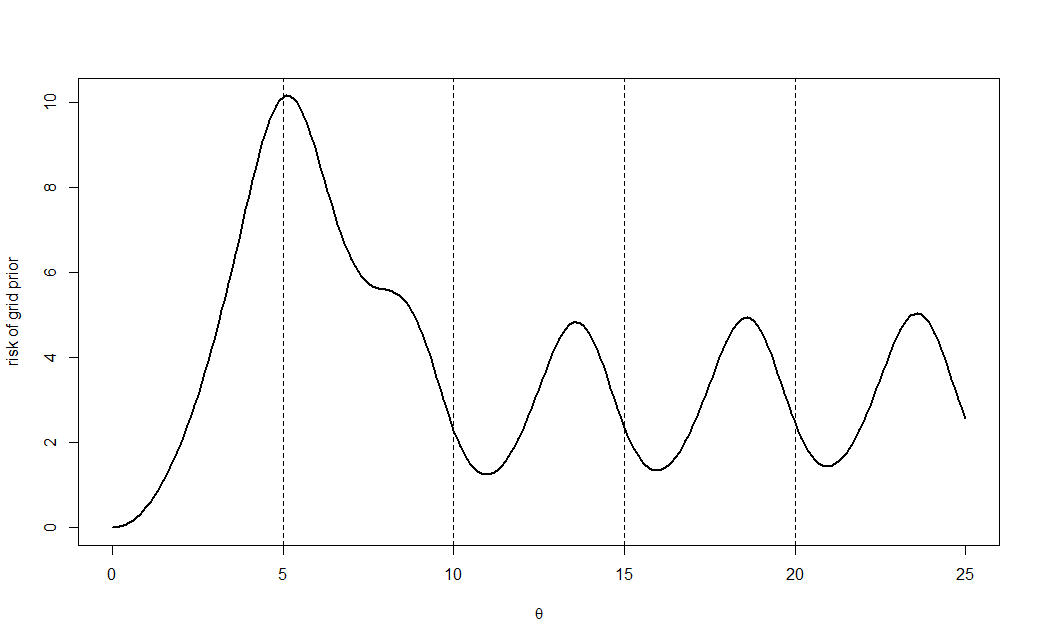}
	\caption{\small \sl Plot of the risk of $\phat_{\g}$ as the parameter $\theta$ varies over $\mathbb{R}^+$. Here, $\lambda=5$, $r=1$ and $\eta=\exp(-\lambda^2)$. The risk is well controlled below the asymptotic theory benchmark minimax value of $\lambda^2/(2r)=12.5$. The vertical line denotes multiples of $\lambda$ along the x-axis.} \label{fig-1}
\end{figure}

\begin{figure}\label{fig-2} 
	\includegraphics[width=\textwidth]{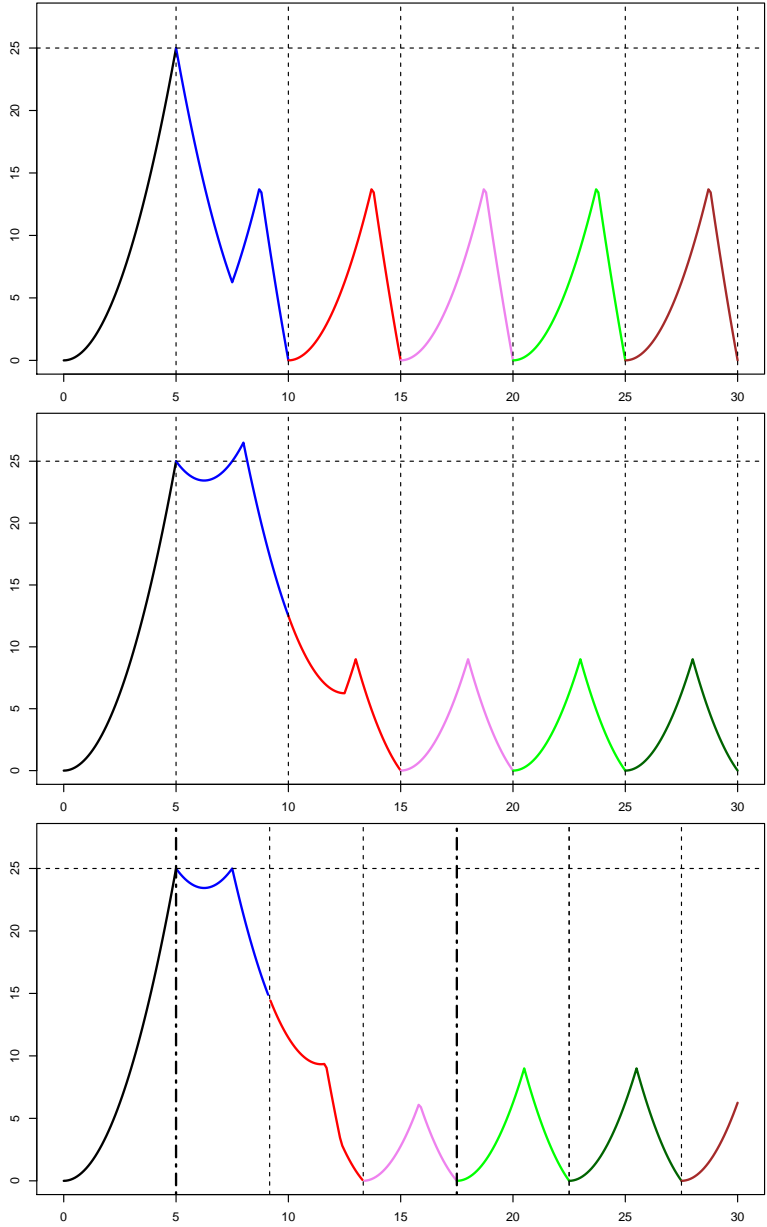}
			\caption{\small \sl
				For $r=1$, the top plot displays the dominant risk $\sigma(l,\omega)$ of $\phat_{\g}$ in the y-axis as $|\theta|=\lambda (\alpha_l +\omega)$ varies in the x-axis. 
				The next two plots show the dominant risk of $\phat_{\g}$ and $\phat_{\bg}$ respectively for $r=0.25$. For each plot, the benchmark minimax value  (adjusted) of  $\lambda^2$ is represented by the dotted horizontal line. The vertical lines partition $|\theta|$ into intervals corresponding to the different values $l$ (based on the co-ordinate system for $\theta$ used in Section~\ref{sec-2}; the purple lines denotes the boundaries for differential spacings) and different colors are used to display the dominant risk in these partitions. Here, $\lambda=5$ for all the 3 plots.  }\label{fig-4}
\end{figure}

To understand the differences in the risk properties between the grid and bi-grid prior as $r$ varies, we now concentrate on the behavior of the risk components theoretically analyzed in Section~\ref{sec-2}. Following the proofs of theorems~\ref{th:D-minimax} and \ref{th:M-minimax},  we concentrate on the dominant component of the risk $\sigma(l,\omega)$ as defined in \eqref{eq:rhobd}  where,  $\theta= \lambda (\alpha_l +\omega)$.
Figure~\ref{fig-2} plots $\sigma(l,\omega)$ for $\phat_{\g}$ and $\phat_{\bg}$ as $l,\omega$ vary. Note, that for the representation of $\theta$ is different for $\pi_{\g}$ and $\pi_{\bg}$ as the $\alpha_l$ and $\beta_l$ are different. As such, for $\phat_{\g}$,  $\omega$ always varies over $[0,1]$ where as for the bi-grid prior $\omega$ only varies over $[0,b]$ for some initial $l$ denominations.  In figure~\ref{fig-2}, we have 3 plots. The top plot shows the dominant risk for $\phat_{\g}$ when $r=1$. The dominant risk here is always contained below the minimax value (shown in dotted horizontal line). Here, $ d^+(l,\omega)$ starts being positive  roughly at $l=1$ and $\omega=0.5$ onwards; and so, the plot of $\sigma(l,\omega)$ for $\phat_{\g}$ has 3 piecewise quadratics for $\theta \in [\lambda,2\lambda]$ but  2 piecewise quadratics for the following intervals.
The next plot shows the dominant risk for $\phat_{\g}$ when $r=0.25$. 
In this case, $ d^+(l,\omega)$ starts being positive  roughly at $l=2$ and $\omega=0.5$ onwards; and so, the plot of $\sigma(l,\omega)$ for $\phat_{\g}$ has 3 piecewise quadratics for $\theta$ in the interval $[2\lambda,3\lambda]$.
Here, the dominant risk is not always contained below the minimax value.  The bottom plot shows  the dominant risk for $\phat_{\bg}$ when $r=0.25$ and it is always contained below the minimax value. Note, that we have used different colors to display the dominant risk  for different values of $l$ and have used vertical lines to partition the $\theta$ values belonging to different values of $l$.  These partitions are same for the top two plots as they both involve $\phat_{\g}$ but is different for $\phat_{\bg}$ in the bottom plot. 

\begin{figure} \begin{center}
		\includegraphics[width=0.9\textwidth]{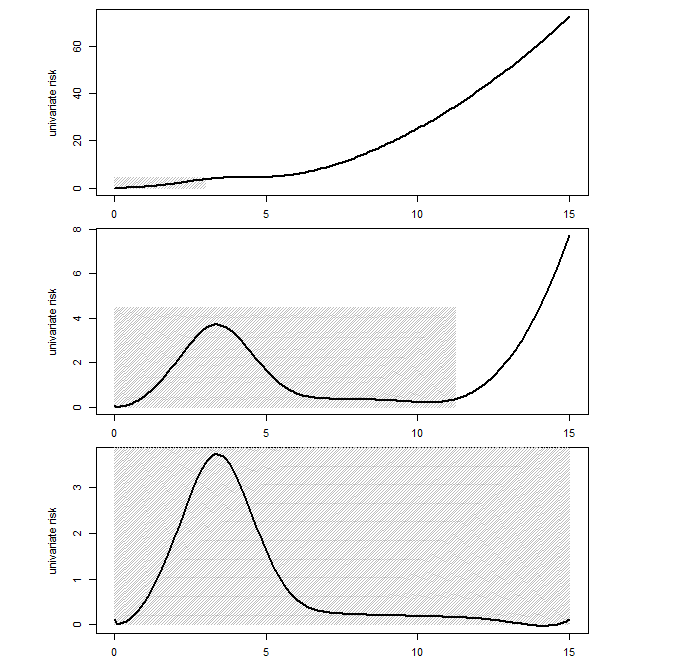}
		\caption{\small \sl From top to bottom we have the plots of the risk of $\phat_S[l]$ as $|\theta|$ varies over the interval $[0,t]$ for $l=\lambda, 0.75t, t$ respectively. Here, $r=1$, $\lambda=3$ and $t=5\lambda$. The shaded rectangle in each plot has base on the support of the prior $[0,l]$ and has height equaling the minimax benchmark value of $\lambda^2/(2r)=4.5$.}\label{fig-big-merge-1}
	\end{center}\label{fig-7}
\end{figure}

\newpage
Figure~\ref{fig-7} displays the risk plots for Spike and Slab prior based \textsf{pdes} $\phat_S[l]$. For three different choices of $l$ the risk plots were tracked the interval $[0,t]$. It was seen that the risk is controlled below the benchmark value of $\lambda^2/(2r)=4.5$ when $l=t$. In that case, the risk function initially increases at a quadratic rate and peaks near $\lambda$ and thereafter decreases at a rapid rate.

Table~\ref{numerical.risk.2} shows the locations of the maxima of the \textsf{pdes} discussed in Section~\ref{sec:simu}. 

\newpage
\begin{table}[!h]
	\centering
		\begin{tabular}{|c|c|c|c|c|c|c|}
			\hline
			\textbf{Sparsity} & \textbf{r} & \textbf{Theory} & \textbf{H-Plugin} & \textbf{C-Thresh} & \textbf{Bi-Grid} & \textbf{SS} \bigstrut\\
		 \hline
		 & 1     & 1.52  & 2.46  & 2.52  & 2.47  & 2.38 \bigstrut\\
		 \cline{2-7}          & 0.5   & 1.24  & 2.41  & 2.49  & 2.17  & 2.10 \bigstrut\\
		 \cline{2-7}    0.1   & 0.25  & 0.96  & 2.37  & 2.46  & 1.80  & 1.78 \bigstrut\\
		 \cline{2-7}          & 0.1   & 0.65  & 2.33  & 2.74  & 1.36  & 1.31 \bigstrut\\
		 \hline
		 & 1     & 2.63  & 3.39  & 3.06  & 3.12  & 3.13 \bigstrut\\
		 \cline{2-7}          & 0.5   & 2.15  & 3.39  & 2.84  & 2.84  & 2.73 \bigstrut\\
		 \cline{2-7}    0.001 & 0.25  & 1.66  & 3.39  & 2.26  & 2.28  & 2.30 \bigstrut\\
		 \cline{2-7}          & 0.1   & 1.12  & 3.39  & 1.68  & 1.69  & 1.70 \bigstrut\\
		 \hline
		 & 1     & 4.80  & 5.93  & 4.94  & 4.91  & 5.05 \bigstrut\\
		 \cline{2-7}    1E-10 & 0.5   & 3.92  & 5.93  & 4.39  & 4.37  & 4.37 \bigstrut\\
		 \cline{2-7}          & 0.25  & 3.03  & 5.93  & 3.43  & 3.89  & 3.58 \bigstrut\\
		 \cline{2-7}          & 0.1   & 2.05  & 5.92  & 2.52  & 2.54  & 2.56 \bigstrut\\
		 \hline
		\end{tabular}	
\caption{\small \sl Numerical evaluation of the location of maxima of the risk plots over $[-l,l]$ for the different univariate predictive densities as the degree of sparsity $(\eta)$ and predictive difficulty $r$ varies. Here, we have chosen $l=5\lambda$, where $\lambda$ is defined in \eqref{eq:lambdadefs}. In `Asymp-Theory' column we report the theoretically obtained first order asymptotic maxima.} 
\label{numerical.risk.2}
\end{table}%


\end{document}